\documentclass[12pt]{amsart}

\usepackage[left=1in,top=1in,right=1in,bottom=1in]{geometry} % sets all margins to 1inch

\input xy
\xyoption{all}
\usepackage{cite} % for bibtex
\usepackage{pstricks}
\usepackage{fancyhdr}
\usepackage{mathrsfs}
\usepackage{epsfig}
\usepackage{amsmath} 
\usepackage{amssymb} 
\usepackage{latexsym}
\usepackage[all, knot]{xy}
\usepackage{verbatim} % for \begin{comment}
\usepackage{tikz} 
\usepackage{wrapfig}
\newtheorem{theorem}{Theorem} %[section]
\newtheorem{lemma}[theorem]{Lemma} 
\newtheorem{proposition}[theorem]{Proposition} %[section]
\newtheorem{corollary}[theorem]{Corollary} %[section]
\theoremstyle{definition}  %[section]
\newtheorem{example}[theorem]{Example} %[section]
\newtheorem{remark}[theorem]{Remark}

\theoremstyle{definition} 
\theoremstyle{definition} 
\theoremstyle{definition} 
\DeclareMathOperator*\argmin{argmin}
\DeclareMathOperator*\argmax{argmax}

\newcommand\C {\mathbb{C}}
\newcommand\R {\mathbb{R}}

\newcommand\PP {\mathbb{P}}
\newcommand\Q {\mathbb{Q}}

\newcommand{\bigO}O

\renewcommand\C{\mathcal{C}}
\renewcommand\L{\mathcal{L}}

\author{Jes\'us A.~De Loera, Bernd Sturmfels, and Cynthia Vinzant}
\email{deloera@math.ucdavis.edu, bernd@math.berkeley.edu, vinzant@umich.edu}
 \address{University of California at Berkeley and Davis; University of Michigan at Ann Arbor}

\title{The Central Curve in Linear Programming}

\begin{document}
\begin{abstract}
The central curve of a linear program is an algebraic curve specified by 
linear and quadratic constraints arising from complementary slackness. 
It is the union of the various central paths for minimizing or maximizing the cost function 
over any region in the associated hyperplane arrangement.
We determine the degree, arithmetic genus and defining prime ideal
of the central curve, thereby answering a question of
Bayer and Lagarias. These invariants, along with the degree of the Gauss image of
the curve, are expressed in terms of the matroid of the input matrix. Extending work of 
Dedieu, Malajovich and Shub, this yields an instance-specific
bound on the total curvature of the central path, a quantity relevant for interior point methods.
The global geometry of central curves is studied in detail.
\end{abstract}
\maketitle

\section{Introduction}
We consider the standard linear programming problem in its
 {\em primal} and {\em dual} formulation:
\begin{equation}
\label{initial problem}
{\rm Maximize} \,\,\,\, \mathbf{c}^T \mathbf{x} \,\,\,\,
	 \textnormal{subject to} \,\,\,
	  A\mathbf{x} = \mathbf{b} \,\,\, {\rm and} \,\,\,
	\mathbf{x} \geq 0;
\end{equation}
\begin{equation}
\label{dual problem slack}
{\rm Minimize} \,\,\,\, \mathbf{b}^T \mathbf{y} \,\,\,\,
         \textnormal{subject to} \,\,\,
          A^T\mathbf{y} - {\bf s}= \mathbf{c} \,\,\, {\rm and} \,\,\, {\bf s} \geq 0.
\end{equation}
Here $A$ is a fixed matrix of rank $d$ having $n$ columns.
The vectors ${\bf c} \in \R^n$ and ${\bf b} \in {\rm image}(A)$ may vary. In most of our
results we assume that ${\bf b}$ and ${\bf c}$ are generic.
This implies  that both the primal optimal solution and the dual optimal solution are  unique.
We also assume that the problem ({\ref{initial problem}) is bounded
and both problems  (\ref{initial problem})  and (\ref{dual problem slack})
are strictly feasible.

Before describing our contributions, we review some basics
   from the theory of linear programming \cite{RTV, vanderbei}.
The (primal) \emph{logarithmic barrier function} for (\ref{initial problem}) is defined as
$$ f_{\lambda}(\mathbf{x}) \,\,\,:= \,\,\, \mathbf{c}^T \mathbf{x} \,+\,
\lambda \sum_{i=1}^n \log x_i, $$ where $\lambda>0$ is a real parameter.  This
specifies a family of optimization problems:
\begin{equation} \label{deform}
{\rm Maximize} \,\,\,\, f_\lambda(\mathbf{x}) \,\,\,\,
	 \textnormal{subject to} \,\,\,
	  A\mathbf{x} = \mathbf{b} \,\,\, {\rm and} \,\,\,
	\mathbf{x} \geq 0 	.
\end{equation}
Since the function $f_\lambda$ is strictly concave, it attains a unique maximum $ \mathbf{x}^{*}(\lambda)$ in the 
interior of the feasible polytope $P = \{ {\bf x} \in \R_{\geq 0}^n : A {\bf x} = {\bf b} \}$. Note that 
$f_\lambda(\mathbf{x})$ tends to $-\infty$ when $\mathbf{x}$ approaches the boundary of $P$.
The {\em primal central path}  is the curve $\{ \mathbf{x}^{*}(\lambda) \,|\, \lambda > 0\}  $ inside the polytope~$P$.
There is an analogous logarithmic barrier function for the dual problem (\ref{dual problem slack})
and a corresponding {\em dual central path}.   The central path connects the optimal solution
of the linear program in question with its \emph{analytic center}. This is the optimal point of $f_\infty$ or equivalently $\argmax_P( \sum_{i=1}^n \log x_i)$.
The central path is homeomorphic to a line segment. 

The {\em complementary slackness} condition says that the pair of optimal solutions, 
to the primal linear program (\ref{initial problem}) and to the dual linear program (\ref{dual problem slack}), 
are characterized by 
\begin{equation}
\label{eq:primaldual}
 A {\bf x} = {\bf b} \,, \,\, A^T {\bf y} - {\bf s}  = {\bf c}\,,\,\, {\bf x} \geq 0 \, ,\,\,
{\bf s} \geq 0 , \,\,\, \text{and} \,\,\,
x_i  s_i = 0 \,\,\, \text{for} \,\,\, i = 1,2,\ldots, n.
\end{equation}
The central path converges to the solution of this system of  equations and inequalities:

\begin{theorem}[cf.~\cite{vanderbei}]
\label{introthm}
If $A$ has $d$ rows, then
for all $\lambda > 0$, the system of polynomial equations 
\begin{equation}
\label{eq:central1}
 A {\bf x} = {\bf b} \,, \,\, A^T {\bf y} - {\bf s}  = {\bf c},\,\, \, \text{and} \,\,\,
x_i  s_i = \lambda \,\,\, \text{for} \,\,\, i = 1,2,\ldots, n,
\end{equation}
 has a unique real solution 
$(\mathbf{x}^*(\lambda), \mathbf{y}^*(\lambda), \mathbf{s}^*(\lambda))$ with the properties
$\mathbf{x}^*(\lambda) > 0$ and $\mathbf{s}^*(\lambda) > 0$. The point
$\mathbf{x}^*(\lambda)$ is the optimal solution of \eqref{deform}. 
The limit point $(\mathbf{x}^*(0), \mathbf{y}^*(0), \mathbf{s}^*(0))$
of these solutions for $\lambda \rightarrow 0$ is the unique solution 
of the complementary slackness constraints \eqref{eq:primaldual}. 
\end{theorem}

\begin{center}
\begin{figure}
\begin{center}
\includegraphics[width=7cm]{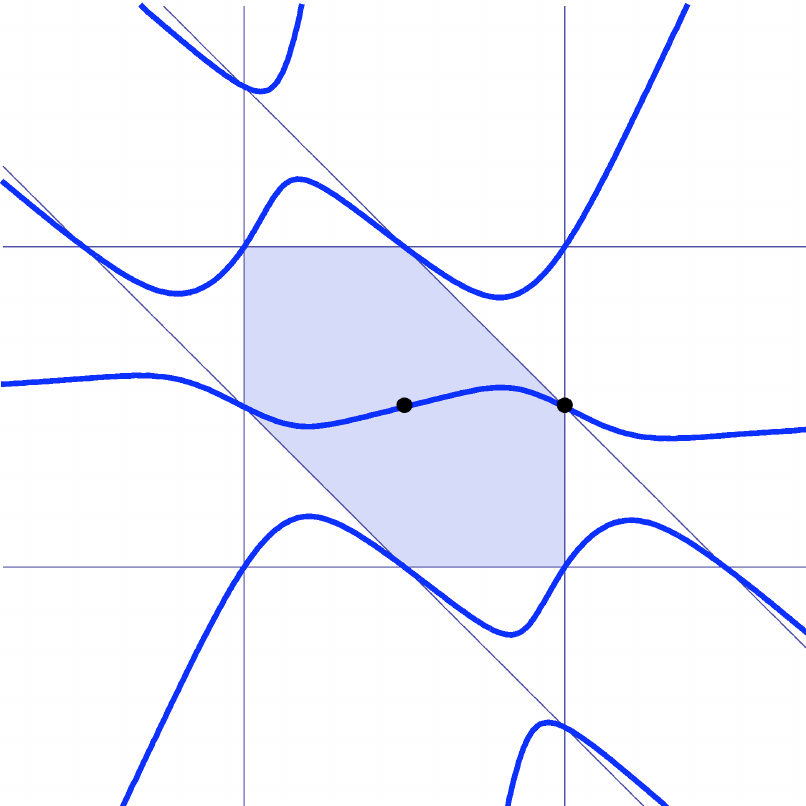}\hspace{1cm}
\includegraphics[width=7cm]{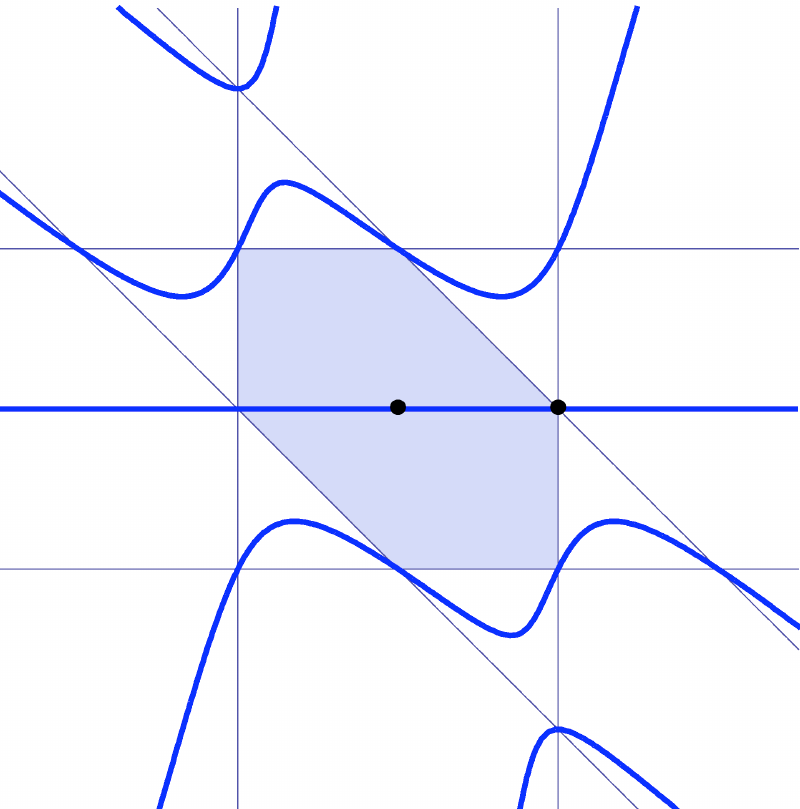}
\end{center}
\caption{The central curve of six lines for two choices of the cost function}
\label{fig:hexagon}
\end{figure}
\end{center}

Our object of study in this paper is the set of {\em all} solutions
of the equations (\ref{eq:central1}), not just those whose coordinates
are real and positive. For general ${\bf b}$ and ${\bf c}$, 
this set is the following irreducible algebraic curve.
The \emph{central curve} is the Zariski closure of
the central path in $({\bf x}, {\bf y}, {\bf s})$-space, that is, it is
 the smallest algebraic variety in $\R^{2n+d}$ that contains
the central path. The {\em primal central curve} in
$\R^n$ is obtained by projecting the central curve into 
${\bf x}$-space. We can similarly define the
{\em dual central curve} by projecting into ${\bf y}$-space or into ${\bf s}$-space.

\begin{example} \label{transport23}
 Figure~\ref{fig:hexagon} depicts the primal central curve for a small
 {\em transportation problem}. Here
 $A $ is the $5 \times 6$ node-edge matrix of the complete bipartite graph $K_{2,3}$, as shown below:\begin{align*}
\begin{tikzpicture}[scale=.8] 
\pgfsetarrowsstart{latex} 
\path (-1,1) node[draw, shape=circle] (v1) {$v_1$};
\path (-1,-1) node[draw, shape=circle] (v2) {$v_2$};
\path (1,2) node[draw, shape=circle] (v3) {$v_3$};
\path (1,0) node[draw, shape=circle] (v4) {$v_4$};
\path (1,-2) node[draw, shape=circle] (v5) {$v_5$};
\draw [-] (v1) -- (v3)
(v2) -- (v4) 
(v2) -- (v3)
(v2) -- (v5) 
(v1) -- (v4)
(v1) -- (v5);
\draw[very thick] [<->] (3,0) -- (5,0);
\path (8,1) node (row1) {$v_1$ \, 1  \,   1 \, 1 \, 0 \, 0 \, 0 };
\path (8,.5) node (row2) {$v_2$ \, 0 \,   0 \, 0  \, 1 \, 1 \, 1 };
\path (8,0) node (row3) {$v_3$ \, 1  \,   0 \, 0 \, 1 \, 0 \, 0 };
\path (8,-.5) node (row4) {$v_4$ \, 0  \,   1 \, 0 \, 0 \, 1 \, 0 };
\path (8,-1) node (row5) {$v_5$ \, 0  \,   0 \, 1 \, 0 \, 0 \, 1 };
\draw[thick, -] (6.5, -1.25) .. controls (6.4, -.5) and (6.4, .5) .. (6.5, 1.25);
\draw[thick, -] (10.4, -1.25) .. controls (10.5, -.5) and (10.5, .5) .. (10.4, 1.25);
\end{tikzpicture} 
\end{align*}
Here $n= 6$ and $d=4$ because $A$ has rank $4$.
We return to this example in Section 4. \hfill $\diamond$
\end{example}

As seen in Figure~\ref{fig:hexagon}, and proved in Theorem \ref{thm:global},
the primal central curve contains the central paths of every polytope
in the arrangement in $\{A {\bf x}={\bf b}\}$ defined 
by the coordinate hyperplanes $\{x_i=0\}$ for the cost functions ${\bf c}$ and $-{\bf c}$. The union over all central curves,
as the right hand side ${\bf b}$ varies, is an algebraic variety
of dimension $d+1$, called the {\em central sheet},
which will play an important role. Our analysis will be based  on results
of Terao \cite{Ter} and Proudfoot--Speyer \cite{PS} on
 algebras generated by reciprocals of linear forms; see also Berget \cite{berget}.
Matroid theory will be our language for working with these algebras and their~ideals.

The algebro-geometric study of central paths was pioneered by
Bayer and Lagarias \cite{BL1, BL2}. Their 1989 articles are part of the 
early history of interior point methods. They observed (on pages 569-571 
of \cite{BL2}) that the central path defines an irreducible algebraic curve 
in ${\bf x}$-space or ${\bf y}$-space, and they identified a complete 
intersection that has the central curve as an irreducible component. The 
last sentence of \cite[\S 11]{BL2} states the open problem of
identifying polynomials that  cut out the central curve, without any extraneous components.
It is worth stressing that one easily finds polynomials that vanish on the central curve from
the gradient optimality conditions on the barrier function. Those polynomials vanish on high-dimensional
components, other than the central curve.
These extra components are contained in the coordinate hyperplanes, and the challenge is to
remove them in our algebraic description.

In numerical optimization, the optimal solution to  \eqref{initial problem} is 
found by following a piecewise-linear approximation to
the central path. Different strategies for generating the
step-by-step moves correspond to different interior point methods. 
One way to estimate the number of Newton steps needed to reach the optimal solution
  is to bound  the \emph{total curvature} of the central
  path. This has been investigated by many authors 
  (see e.g.~\cite{DMS, MonteiroTsuchiya, Sonnevendetal,VavasisYe,ZhaoStoer}),
  the idea being that  curves with small curvature are easier to approximate with line segments.
  The algebraic results in this paper
    contribute to the understanding of the total curvature.

Here is a list of our results. Precise statements are given in each section.

\begin{itemize}

\item 
Section  \ref{planarcurves} analyzes central curves in the plane, 
with emphasis on the dual formulation $(d=2)$. We show that our
 curves are {\em Vinnikov curves} \cite{Vinnikov} of degree $\leq n-1$, obtained 
from an arrangement of $n$ lines by taking a {\em Renegar derivative} \cite{Ren2}. The  total
curvature of a plane curve can be bounded in terms of its number of real inflection points.
We  derive a new bound from a classical formula due to Felix~Klein~\cite{Klein}.

\item  All our formulas and bounds in Sections \ref{centralpathideal}, \ref{gaussc}, \ref{dual_and_average}, and \ref{globalgeo}   are expressed in the language of matroid theory. A particularly important role is played by matroid
 invariants, such as the {\em Tutte polynomial}, that are associated with the matrix $A$. In
 Section \ref{matroidstuff} we review the required background
  from matroid theory and geometric combinatorics. 

\item In Section \ref{centralpathideal}  we present a complete solution to
the Bayer-Lagarias problem.  Under the assumption that
${\bf b}$ and ${\bf c}$ are general, while $A$ is fixed  and possibly special,
we determine the prime ideal of all polynomials that vanish on
the primal central curve.  We express the degree of this curve as a
matroid invariant.  This yields the tight upper bound
$\binom{n-1}{d}$ for the degree. For instance, the curves in
Figure~\ref{fig:hexagon} have degree five.
 We also determine the Hilbert series and arithmetic genus of our curve in~$\PP^{n}$.

\item Section~\ref{gaussc} develops our approach to estimating the total curvature of the central curve.
  Dedieu, Malajovich and Shub \cite{DMS} noted 
 that the total curvature of any curve $\mathcal{C}$  coincides with the 
 arc length of the image of $\mathcal{C}$ under the \emph{Gauss map}.
 Hence any bound on the degree of the {\em Gauss curve} translates into
 a bound on the total curvature. Our main result in Section~\ref{gaussc}
 is a very precise bound, in terms of matroid invariants,
  for the degree of the Gauss curve arising from any linear program. 

\item While Sections  \ref{centralpathideal} and \ref{gaussc}
focused on primal linear programs.
Section \ref{dual_and_average} revisits our results on the degree and 
curvature, and it translates them to the dual formulation.
Theorem \ref{avgCurv} characterizes the
average curvature over the bounded feasibility regions.

\item  Section \ref{globalgeo} furnishes an entirely symmetric
description of the primal-dual central curve inside a product of
two projective spaces.  This leads to a range of results on
the global geometry of our curves. In particular, we explain
 how the central curve passes through all
vertices of the hyperplane arrangement and through all
the  analytic centers.
 \end{itemize}

\smallskip

What got us started on this project was our desire to understand the ``snakes'' 
of Deza, Terlaky and Zinchenko \cite{DTZ08}. We close the introduction by presenting their curve for $n=6$.

\begin{example}  \label{ex:snake}
Let $n=6$, $d=2$ and fix the following matrix, right hand side and cost vector:
$$ A = \begin{pmatrix} 0 &    -1 &     1 &    -1  &    1 &      -1  \\
-1 & \frac{1}{10} & \frac{1}{3} & \frac{100}{11} & \frac{1000}{11}  & \frac{10000}{11} 
\end{pmatrix} \,, \quad
{\bf b} \,=\,\begin{pmatrix} 0 \\ 1 \end{pmatrix} , $$
$${\bf c}^T =  \begin{pmatrix} -1 &  -\frac{1}{2} & - \frac{1}{3} & 
-\frac{449989}{990000} & 
-\frac{359989}{792000} & 
-\frac{299989}{660000} \end{pmatrix} . $$
The resulting linear program, in its dual formulation (\ref{dual problem slack}), is precisely the
instance in \cite[Figure 2, page 218]{DTZ08}. We redrew the central curve in Figure \ref{fig:dtz6}.
The hexagon $P_{6,2}^*$ shown there equals $\,\{{\bf y} \in \R^2 \,:\,  A^T{\bf y} 
 \geq {\bf c} \}$.
 The analytic center of $P_{6,2}^*$ is  a point  with  approximate coordinates
$\, {\bf y} =(-0.027978 ..., 0.778637...)$. It has
 algebraic degree $10$ over $\Q$, which indicates the level of difficulty to write exact coordinates.
  The optimal solution is the vertex with rational coordinates  $\,{\bf y}= 
  (y_1,y_2) = (-\frac{599700011}{1800660000},  -\frac{519989}{600220000} )
 =  (-0.033304..., -0.00086...).$

Following \cite{DTZ08},
%the authors used a sample of points to estimate the value
%of the total curvature of the central path. W
we sampled many points along the central
path, and we found that
%obtained the following more accurate figure: 
%The arc length of the central path is at least $\,1.1149315669943... $ and 
the total  curvature of the central path equals $13.375481417...$ .  This measurement concerns only the 
part of the central curve that goes from the analytic center to the optimum.  
%Using the
%techniques of this paper we can provide an upper bound to the total curvature (see Section \ref{gaussc}).
Our algebraic recipe (\ref{Jdet}) for computing the central curve leads to the following polynomial:
{\tiny
$ \left( y_{{2}}-1 \right)  \bigl( 2760518880000000000000000\,y_2^{4}+22783991895360000000000000\,y_2^3 
y_1-1559398946696532000000000\,y_2^3 + \\ 1688399343321073200000000 y_1 y_2^2+87717009913470910818000 y_2^2-3511691013758400000000000 y_1^2 y_2^2-324621326759441931317 y_2 \\ +11183216292449806548000 \,y_1 y_2
+2558474824415400000000\, y_1^2 y_2 - 51358431801600000000000\, y_1^3 y_2+
6337035495096700140\,y_1 \\ + 77623920000000000000\, y_1^4-13856351760343620000\,y_1^2 +291589604847546655-
38575873512000000000\, y_1^3\bigr).
$
}\smallskip \\
This polynomial of degree five has a linear factor $y_2-1$ because the vector ${\bf b}$
that specifies the objective function in this dual formulation
 is parallel to the first column of $A$. Thus the central curve in Figure \ref{fig:dtz6} has degree four,
and its defining irreducible polynomial is the second factor.
When the cost vector ${\bf b}$ is replaced by a vector that is not parallel to
a column of $A$ then the output of the same calculation
(to be explained in Section 4) is an irreducible polynomial of degree five.
In other words, for almost all ${\bf b}$, the central curve is a quintic.

\begin{figure}
\includegraphics[width=6.4cm]{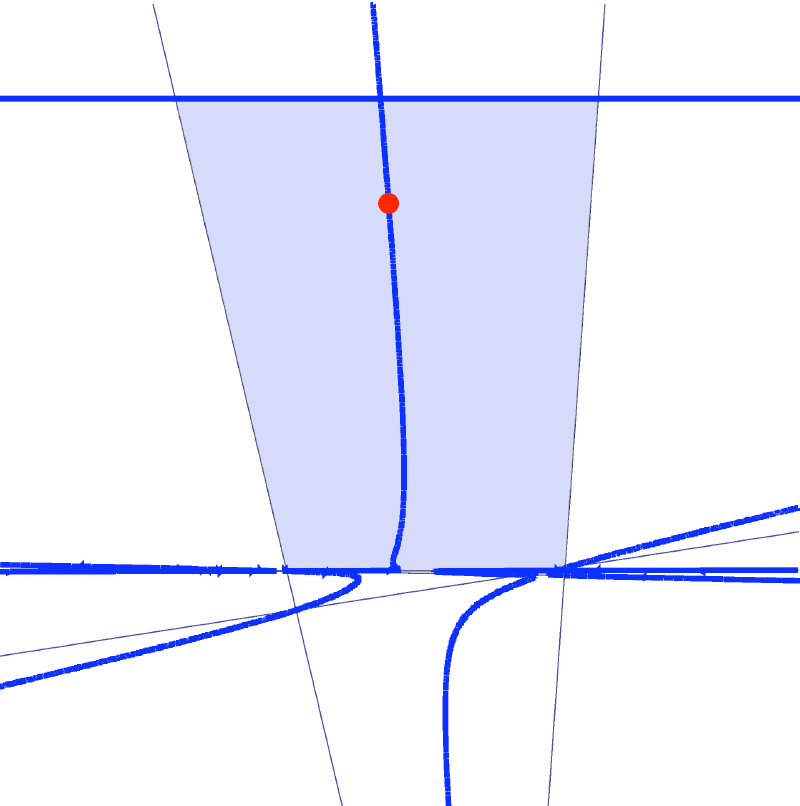}\hspace{1.3cm}
\includegraphics[width=6.4cm]{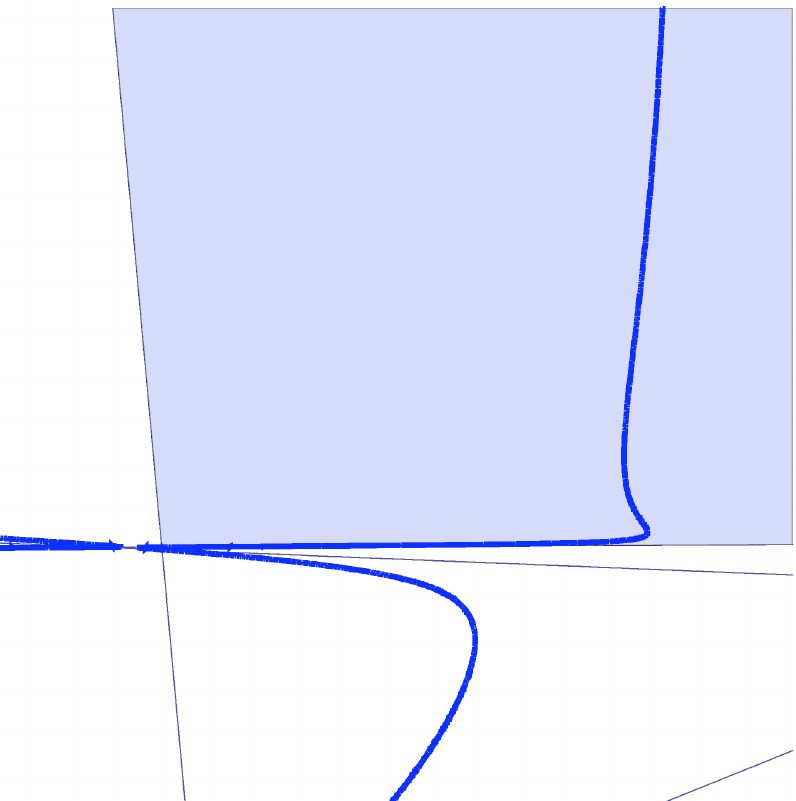}
\caption{The DTZ snake with 6 constraints. 
On the left, a global view of the polygon and its central
 curve with the line $y_2=1$ appearing 
as part of the curve.
On the right a close-up of the central path and its inflection points.}
\label{fig:dtz6}
\end{figure}

While most studies in optimization focus only on just the small portion of
the curve that runs from the analytic center to the optimum, we 
argue here that the algebraic geometry of the entire curve reveals a
more complete and interesting picture.  For generic ${\bf b}$ and ${\bf c}$, the  central curve is a
quintic that passes through all vertices of the line arrangement
defined by the six edges of the polygon. As we shall see, it
passes through the analytic centers of all bounded cells
(Theorem \ref{thm:global})
 and it is topologically a nested set
of ovals (Proposition \ref{prop:hyperbolicpolygon}). \hfill $\diamond$
\end{example}

\section{Plane Curves}
\label{planarcurves}

When the central curve lives in a plane, the curve is cut out by a single polynomial equation.
This occurs for the dual curve when $d=2$ and the primal curve when $n=d-2$. We now focus 
on the dual curve ($d=2$). This serves as a warm-up to the full derivation of all 
equations in Section~\ref{centralpathideal}. In this section we derive the equations of the central 
curve from first principles, we show that these curves are hyperbolic and Renegar derivatives of 
products of lines, and we use this structure to bound the average total curvature of the curve. 

Let $A=(a_{ij})$ be a fixed $2\times n$ matrix of rank $2$, and consider
arbitrary vectors ${\bf b} = (b_1,b_2)^T \in \R^2$ and ${\bf c} = (c_1,\ldots,c_n)^T \in \R^n$.
Here the ${\bf y}$-space is the plane $\R^2$ with coordinates ${\bf y}=(y_1,y_2)$.
The central curve is the Zariski closure in this plane of the parametrized path
\[{\bf y}^*(\lambda) \;= \;\argmin_{\{{\bf y}\;:\; A^T{\bf y}\geq {\bf c}\}}\;\; b_1y_1+b_2y_2 - \lambda\sum_{i=1}^n \log(a_{1i}y_1 +a_{2i}y_2-c_i). \]
The conditions for optimality are obtained by
setting the first partial derivatives to zero:
$$ 0\; =\; b_1 -\lambda  \sum_{i=1}^n \frac{a_{1i}}{a_{1i}y_1 +a_{2i}y_2-c_i} \;\;\;\;\text{  and  }\;\;\;\; 0 \;=\; b_2 - \lambda  \sum_{i=1}^n \frac{a_{2i}}{a_{1i}y_1 +a_{2i}y_2-c_i}.$$ 
Multiplying these equations by $b_2/\lambda$ or $b_1/\lambda$ gives 
\begin{equation}
\label{beforecurvaplana}
\frac{b_1b_2}{\lambda}\; \;= \;\;  \sum_{i=1}^n \frac{b_2a_{1i}}{a_{1i}y_1 +a_{2i}y_2-c_i}
\; \;= \;\;  \sum_{i=1}^n \frac{b_1a_{2i}}{a_{1i}y_1 +a_{2i}y_2-c_i}.
\end{equation}
This eliminates the parameter $\lambda$ and we are left with the equation on the right.
By clearing denominators, we get a single polynomial $C$ that vanishes on 
 the central curve in ${\bf y}$-space: 
\begin{equation}
\label{curvaplana}
 C({\bf y})\;\;= \;\; \sum_{i\in \mathcal{I}} (b_1a_{2i}-b_2a_{1i})
 \prod_{j\in \mathcal{I}\backslash\{i\}} (a_{1j}y_1 +a_{2j}y_2-c_j),
 \end{equation}
 where $\mathcal{I} = \{i\;:\; b_1a_{2i}-b_2a_{1i}\neq 0\}$. We see that 
  the degree of $C({\bf y})$ is $|\mathcal{I}|-1$. This  equals $n-1$ for generic ${\bf b}$.
  In our derivation we assumed that $\lambda$ is non-zero but
the resulting equation is valid on the Zariski closure, which includes the important
points with parameter $\lambda=0$. 

We consider the closure $\mathcal{C}$ of the central curve in
the complex projective plane $\PP^2$ with coordinates $[y_0:y_1:y_2]$.  
Thus $\mathcal{C}$ is the complex projective curve defined by
$\,y_0^{|\mathcal{I}|-1}C(\frac{y_1}{y_0},\frac{y_2}{y_0})$.

\begin{proposition}\label{prop:hyperbolicpolygon}
The curve $\mathcal{C}$ is \emph{hyperbolic} with respect to the point $[0:-b_2:b_1]$. This means that every line
 in $\PP^2(\R)$  passing through this special point  meets $\mathcal{C}$ only in real points. \end{proposition}

\begin{proof}
Any line passing through the point  $[0:-b_2:b_1]$ (except the line $\{y_0=0\})$ has the form
$\{b_1y_1+b_2y_2 = b_0y_0\}$ for some $b_0 \in \R$. 
See the left picture in Figure \ref{fig:polygons}. We shall see in
 Remark~\ref{rem:CPtoAC} that,  for any $b_0 \in \R$, the line meets $\mathcal{C}$ in 
 $\deg(\mathcal{C})$ real points.
\end{proof}

\begin{figure}
\includegraphics[width=7cm]{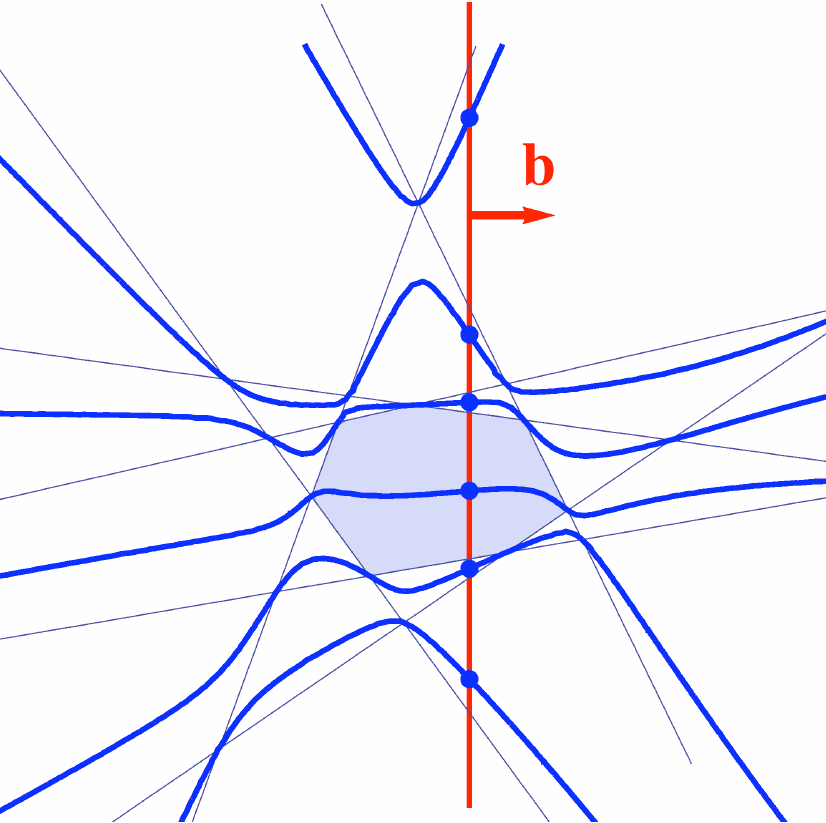} \hspace{1cm}
\includegraphics[width=7cm]{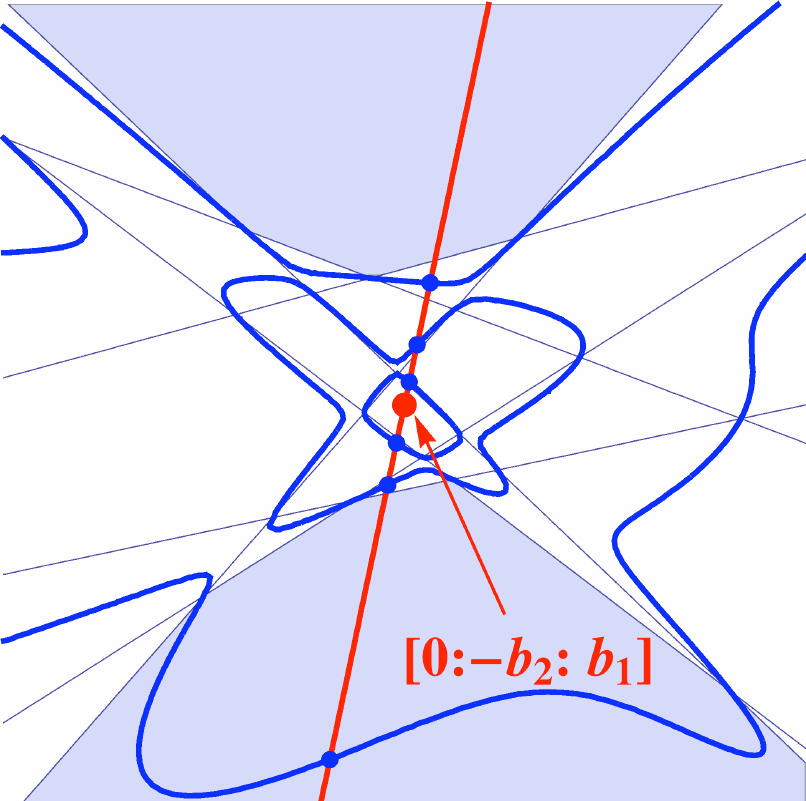}
\caption{The degree-6 central path of a planar 7-gon 
in the affine charts $\{y_0=1\}$ and $\{y_2=1\}$. Every line passing through $[0:-b_2:b_1]$ intersects
the curve in 6 real points, showing the real curve to be 3 completely-nested ovals. }
\label{fig:polygons}
\end{figure}

Hyperbolic curves are also known as {\em Vinnikov curves}, in light of
Vinnikov's seminal work \cite{LPR, Vinnikov} relating them
to semidefinite programming \cite{philipbernd}.
Semidefinite programming has been generalized to hyperbolic programming,
in the work of Renegar \cite{Ren2} and others.
A key construction in hyperbolic programming is the Renegar derivative
which creates a (hyperbolic) polynomial of degree $D-1$
from any (hyperbolic) polynomial of degree $D$.
To be precise, the \textit{Renegar derivative} of a homogeneous
polynomial $f$ with respect to a point~${\bf e}$~is  
\[R_{\bf e}f({\bf y}) \,\,\, = \,\,\, \left( \frac{\partial}{\partial t} f({\bf y}+ t {\bf e}) \right)\biggl|_{t=0}.\]
Renegar derivatives correspond to the {\em polar curves}
 of classical algebraic geometry \cite[\S 1.1]{Dolgachev}.

The Renegar derivative of $f= \prod_{i\in \mathcal{I}}(a_{1i}y_1 +a_{2i}y_2-c_iy_0)$ 
with ${\bf e }= (0,-b_2,b_1)$ is seen to~be
\begin{equation}
\label{eq:renegar}
R_{\bf e}f({\bf y}) \quad = \quad
 \sum_{i\in \mathcal{I}} (b_1a_{2i}-b_2a_{1i})\!\!\prod_{j\in \mathcal{I}\backslash \{i\}}
 \!\! (a_{1j}y_1 +a_{2j}y_2-c_jy_0) \;\;=\;\;C({\bf y}).
 \end{equation}
 In words:~the central curve $\mathcal{C}$ is the Renegar derivative,
taken with respect to the cost function, of the
product of the linear forms that define the convex polygon of feasible points.

The product of linear forms $f = \prod_i (a_{1i}y_1 +a_{2i}y_2-c_iy_0)$
is a hyperbolic polynomial with respect to ${\bf e}$.
 Renegar \cite{Ren2} shows that if $f$ is hyperbolic with respect to ${\bf e}$ then so is $R_{\bf e}f$.
 This yields a second proof  for Proposition \ref{prop:hyperbolicpolygon}. 

Proposition \ref{prop:hyperbolicpolygon} is visualized in Figure~\ref{fig:polygons}.
The picture on the right is obtained from the picture on the left by
a projective transformation. The point at infinity which represents
the cost function is now in the center of the diagram.
In this rendition, the central curve consists of three nested ovals around that point,
highlighting the salient features of a Vinnikov curve.
This beautiful geometry is found not just in the dual picture
but also in the primal picture:

\begin{remark} \label{hyperbolicprimal}
If $d=n-2$ then the primal central curve lies in the plane $\{A{\bf x}={\bf b}\}$. 
The conditions for optimality of (\ref{initial problem}) state that the vector $\nabla(\sum_i\log x_i) = (x_1^{-1}, \hdots, x_n^{-1})$ is in the span of ${\bf c}$ and the rows of $A$. 
The Zariski closure of such vectors is the \textit{central sheet},
to be seen in Section 4. Here, the
central sheet is the hypersurface in $\R^n$ with defining polynomial
\begin{equation}
\label{eq:sec2matrix}
\; \det\begin{pmatrix} 
A_1 \! & \! A_2 \! & \!  \cdots & \! A_n \\
c_1 \! & \!  c_2 \! & \! \cdots \! & \! c_n\\
 x_1^{-1} & x_2^{-1} & \cdots  & x_n^{-1} \\
\end{pmatrix} \cdot \prod_{i\in \mathcal{I}} x_i,
\end{equation}
where $A_i$ is the $i$th column of $A\,$ and $\,\mathcal{I} = \{\, i\;:\;\{\binom{A_j}{{\bf c}_j}\}_{ j\in {[n]\backslash i}}\text{ are linearly independent}\}$.  We see that the degree of this hypersurface is
 $|\mathcal{I}|-1$, so it is  $n-1$ for generic $A$. Intersecting this surface with the plane 
 $\{A {\bf x}={\bf b}\}$ gives the primal central curve, which is hence a curve of degree
 $|\mathcal{I}|-1$.   The corresponding complex projective  curve in 
 $ \PP^2 =\{\,[x_0:{\bf x}] \,|\,A {\bf x}=x_0 {\bf b} \} \subset \PP^{n}$ 
  is hyperbolic with respect to the point $[0: {\bf v}]$ in $\PP^{n}$, 
  where ${\bf v}$ spans the kernel of $\binom{A}{{\bf c}}$. \hfill $\diamond$
\end{remark}

It is of  importance for interior point algorithms to know the exact total curvature, 
formally introduced in equation (\ref{curvaturedef}),
of the central path of a linear program (see \cite{DMS, MonteiroTsuchiya,Sonnevendetal, VavasisYe,ZhaoStoer}).
Deza {\it et al.}~\cite{DTZ08} proved that even for $d=2$ the total curvature grows 
linearly in $n$, and they conjectured that the total curvature is no more than $2\pi n$. 
They named this conjecture the \emph{continuous Hirsch conjecture} because of its 
similarity with the discrete simplex method analogue (see \cite{DTZ}). In Section~\ref{gaussc} we derive general bounds for total curvature, but for plane curves we can exploit an additional geometric
feature, namely, {\em inflection points}.

Benedetti and Ded\`o \cite{Ben} derived a general bound for
the total curvature of a real plane curve in terms of its number of inflection points and its degree. 
We can make this very explicit  for our central path $\{{\bf y}^*(\lambda):\; \lambda \in \mathbb{R}_{\geq 0}\}$.
Its total curvature is bounded above by
\begin{equation}
\label{eq:inflectiobound}
 \text{total curvature of the central path}  \quad \leq \quad \pi\cdot (\text{its number of inflection points} +1). 
 \end{equation}
 To see this, consider the Gauss map $\gamma$ that  takes the curve into the circle $S^1$ 
 by mapping any point on the curve to its unit tangent vector.
  See Section \ref{gaussc} for the general definition.
  The total curvature is the arc length of the image of the Gauss map. 
As $\lambda$ decreases from $\infty$ to 0, the cost function ${\bf b}^T {\bf y^*(\lambda)}$ 
strictly decreases. This implies that, for any point ${\bf y}^*(\lambda)$ on the curve, its image
 under the Gauss map has positive inner product with ${\bf b}$, that is, 
 ${\bf b}^T\gamma({\bf y}^*(\lambda)) \geq 0$. Thus the image of the Gauss map is 
 restricted to a half circle of $S^1$, and it cannot wrap around $S^1$. This shows that the 
 Gauss map can achieve a length of at most $\pi$ before it must ``change direction", 
which happens only at inflection points of the curve. 

It is known that the total number of (complex) inflection points of a plane 
curve of degree $D$ is at most $3D(D-2)$. For real inflection points, there is an even better bound:

\begin{proposition}[A classical result of Felix Klein \cite{Klein}] \label{klein}
\hfill \break The number of real inflection points of a  plane curve of degree $D$ is at most $D(D-2)$.
\end{proposition}

This provides only a quadratic bound for the total curvature of the central path in terms of its degree, 
but it does allow us to improve known bounds for the average total curvature. The \emph{average total curvature} 
of the central curve of a hyperplane arrangement is the average, 
over all bounded regions of the arrangement, of the total curvature of the central curve in that region.
Dedieu {\it et~al.} \cite{DMS} proved that the average total curvature in a simple arrangement 
({\it i.e.}~for a generic matrix $A$) defined by $n$ hyperplanes in dimension $d$ is not greater than $2\pi d$. 
When $d=2$, we can use Proposition~\ref{klein} to improve this bound by a factor of two. See Theorem~\ref{avgCurv} for the 
 extension to general $d$. 

\begin{theorem} The average total curvature of a central path of a generic line arrangement in the plane is at most $2\pi$. 
\end{theorem} 

\begin{proof} The central curve for $n$ general lines
in $\R^2$ has degree $n-1$ and consists of $n-1$ (real affine) connected components.
The argument above and Klein's theorem then show that
\begin{align*}\sum_{i=1}^{n-1}(\text{curvature of the $i$th component})&\,\leq \, \sum_{i=1}^{n-1}\pi(\# \text{inflection points on the $i$th component}+1)\\
&\,\leq \,\,\, \pi (n-1)(n-2). \end{align*}
Our arrangement of $n$ general lines  has
 $\binom{n-1}{2}$ bounded regions. The average total curvature over each of these regions is 
 therefore at most  $\,\pi (n-1)(n-2) / \binom{n-1}{2} \,=\, 2\pi$.  
\end{proof}

To bound the curvature of just the central path, we need to bound the number of inflection points appearing 
on that piece of the central curve.  To address this issue,
we posed the following problem in the manuscript version of this article:
{\em 
 What is the largest number of inflection points on a single oval of a hyperbolic
 curve of degree $D$ in the real plane? In particular, is this number {\em linear} in the degree $D$?
 }
These questions have since been answered by Erwan Brugall\'e and 
Luc\'ia L\'opez de Medrano, using an extension of their techniques in \cite{BruMed}. 
They constructed a Vinnikov curve of even degree $D$ which has 
the maximal number $D(D-2)$ of inflection points
and all of these inflection points lie on the outermost oval.
It would be very interesting to see whether their approach
can be applied to improve the
DTZ snakes of Example \ref{ex:snake}
and lead to   new lower bounds for the total curvature
of planar central paths.

\section{Concepts from Matroid Theory}
\label{matroidstuff}

We have seen in the previous section that the geometry of a central curve
in the plane is intimately connected to that of the underlying arrangement of
constraint lines. For instance, the degree of the central curve, $|\mathcal{I}|-1$, is one less
than the number of constraints not parallel to the cost function.
The systematic study of this kind of combinatorial information,
encoded in a geometric configuration of vectors or hyperplanes, is the subject
of {\em matroid theory}.

Matroid theory will be crucial for stating and proving our results in the
rest of this paper. This section offers an exposition of the relevant concepts.
Of course, there is already a well-established connection between matroid theory and 
linear optimization (e.g., as outlined in \cite{lawler} or in
oriented matroid programming \cite{bachemkern}). Our paper sets up yet another connection.
The material that follows is well-known in algebraic combinatorics, 
but less so in optimization, so we aim to cover the basic facts. The missing 
details can be found in \cite{bjornercomplex,tuttesurvey}.

We consider an $r$-dimensional  linear subspace $\L$ of the vector space $K^n$ with its fixed standard basis. 
Here $K$ is any field.  Typically, $\L$ will be given to us as the row space of an $r \times n$-matrix.
The kernel of that matrix is denoted by $\L^\perp$. This is  a subspace of dimension $n-r$ in $K^n.$
We write $x_1,\ldots,x_n$ for  the restriction of the standard coordinates on $K^n$ to  $\L$.

The two subspaces $\L$ and $\L^\perp$ specify a dual pair of matroids,
denoted $M(\L)$ and $M(\L^\perp)$, on the set $[n] = \{1,\ldots,n\}$.
The matroid $M(\L)$ has rank $r$ and its dual $M(\L^\perp) = M(\L)^*$
has rank $n-r$. We now define the first matroid $M = M(\L)$
by way of its {\em independent sets}. A subset $I$ of $[n]$ is 
{\em independent} in $M$ if the linear forms in $\{x_i : i \in I\}$ 
are linearly independent  on $\L$. Maximal independent sets are called {\em bases}. These
all have cardinality $r$. A subset $I$ is {\em dependent} if it
is not independent. It is a {\em circuit} if it is minimally dependent.

\begin{example} 
Consider the linear space $\L$ spanned by the rows of the rank $4$ matrix $A$ in 
Example~\ref{transport23}.  Because the first four columns of $A$ are linearly independent, the linear
forms $\{x_1, x_2, x_3, x_4\}$ are linearly independent on $\L$ and $\{1,2,3,4\}$ is an independent set 
of $M(\L)$.  As $\L$ has dimension four, $\{1,2,3,4\}$ is a basis of the matroid $M(\L)$.  
On the other hand, the set of columns $\{1,2,4,5\}$ is linearly dependent but every
proper subset is linearly independent. Hence
$\{1,2,4,5\}$ is a circuit of $M(\L)$. From similar considerations, we find that 
$M(\L)$ has nine circuits, namely $\{i,j,k,l\}$ where $i,j\in\{1,2,3\}$ and $k,l\in \{4,5,6\}$.\hfill $\diamond$
\end{example}

One matroid application of importance for our study of central curves 
is the following formula for number of bounded components of a hyperplane arrangement. 
Let ${\bf u} $ be a generic vector in $\R^n$ and consider the
$(n-r)$-dimensional affine space $\L^{\perp} + {\bf u}$ of $\R^n$.
The equations $x_i = 0$ define $n$ hyperplanes in this affine space.
The  arrangement $\{x_i=0\}_{i \in [n]}$  in $\L^{\perp}+ {\bf u}$
is {\em simple}, which means that no point lies on more than $n-r$ of the $n$ hyperplanes.
The vertices of this hyperplane arrangement are in bijection with the bases of the matroid $M$.
The complements of the hyperplanes are convex polyhedra; they are the {\em regions} of the arrangement.
Each region is either bounded or unbounded, and we are interested in the bounded regions. These bounded
regions are the feasibility regions for the linear programs with various sign restrictions on the variables $x_i$. 
Proposition 6.6.2 in \cite{tuttesurvey}, which is 
based on results of Zaslavsky \cite{zaslavsky1}, equates the number of such regions with the 
absolute value of the {\em M\"obius invariant} $\mu(M)$ of the matroid of $M$:
\begin{equation}
\label{eq:moebius4}
|\mu(M)| \,= \,
\# \,\hbox{bounded regions of the hyperplane arrangement $\{x_i=0\}_{i \in [n]}$  in $\L^{\perp}{+} {\bf u}$.} \!
 \end{equation}
Further below, in Equation (\ref{eq:moebius5}),
this invariant will be expressed in terms of the matroid $M$. 
We refer to $|\mu(M)|$   as the {\em M\"obius number} of the matroid~$M$.  
 
To obtain the M\"obius number and more refined invariants that we will need in Section~\ref{gaussc} 
we  introduce a simplicial complex associated to the matroid $M$ called the {\em broken circuit complex}.
We fix the standard ordering $1 < 2 < \cdots < n$ of $[n]$. A {\em broken circuit} of $M$
is any subset of $[n]$ of the form $C \backslash \{{\rm min}(C)\}$ where
$C$ is a circuit. The \emph{broken circuit complex} of $M$ is
the simplicial complex ${\rm Br}(M)$ whose minimal non-faces are the 
broken circuits. Hence, a subset of $[n]$ is a face of ${\rm Br}(M)$ if it does
 not contain any broken circuit. It is known that ${\rm Br}(M)$ is a shellable 
 simplicial complex of dimension $r-1$ (see Theorem 7.4.3 in \cite{bjornercomplex}).
We can recover the M\"obius number of $M$ as follows.
Let $f_i = f_i({\rm Br}(M))$ denote the number of $i$-dimensional faces of the broken
circuit complex ${\rm Br}(M)$. The corresponding h-vector 
$(h_0,h_1,\ldots,h_{r-1})$ can be read off from any shelling (cf.~
\cite[\S 7.2]{bjornercomplex} and \cite[\S 2]{stanley}). It satisfies
\begin{equation}
\label{eq:FH}
 \sum_{i=0}^{r-1} \frac{f_{i-1} z^i }{(1-z)^i} \quad = \quad
\frac{h_0 + h_1 z + h_2 z^2 + \cdots + h_{r-1} z^{r-1}}{(1-z)^r} .
\end{equation}
The relation between the f-vector and the h-vector holds for any simplicial complex \cite{stanley}. 
The rational function (\ref{eq:FH}) is the {\em Hilbert series} (see \cite{stanley}) of the
{\em Stanley-Reisner ring} of the broken circuit complex ${\rm Br}(M)$. The defining ideal of the
Stanley-Reisner ring is generated by the monomials
$\,\prod_{i \in C \backslash \{{\rm min}(C)\}} x_i\,$ representing broken circuits.
Proudfoot and Speyer \cite{PS} constructed a {\em broken circuit ring},
which is the quotient of $K[x_1,\ldots,x_n]$
modulo a prime ideal whose initial ideal is 
precisely this monomial ideal. Hence (\ref{eq:FH}) is also the Hilbert series of the ring in \cite{PS}.
In particular, the  M\"obius number is the common degree of both rings:
\begin{equation}
\label{eq:moebius5}
|\mu(M)| \quad = \quad h_0 + h_1 + h_2 + \cdots + h_{r-1}.
\end{equation}

\begin{example}[Uniform matroids] \label{ex:uniform}
  If $\L$ is a general $r$-dimensional subspace of $K^n$
then $M=M(\L)$ is  the {\em uniform matroid} $M = U_{r,n}$, whose
bases are all $r$-subsets in $[n]$ and whose circuits are all
$(r+1)$-subsets of $[n]$. The broken circuits of $M$ are then all the $r$-subsets
of $\{2,\hdots, n\}$. The broken circuit complex ${\rm Br}(M)$ is the $(r-1)$-dimensional simplicial complex on $[n]$ whose maximal simplices
are $\{1,j_1,\ldots,j_{r-1}\}$ where $2\leq j_1 < \cdots < j_{r-1} \leq n$.
This shows that  $f_i$ equals $ \binom{n-1}{i}$ for $1\leq i \leq r-1$. Using \eqref{eq:FH}
one finds that $\,h_i =  { n-r+i-1 \choose i }$. 
We can then use \eqref{eq:moebius5} to compute the M\"obius number of the uniform matroid $M = U_{r,n}$:
\begin{equation}
\label{eq:moebius2}
|\mu(U_{r,n})|\quad = \quad \sum_{i=0}^{r-1}  \binom{n-r+i-1}{i} \quad = \quad \binom{n-1}{r-1}.
\end{equation}
This binomial coefficient is an upper bound on $|\mu(M)|$ for any  rank $r$ matroid $M$ on $[n]$.

To understand the geometric interpretation of $|\mu(M)|$, let us identify $\L^{\perp}$ with
$\R^{n-r}$.  Here we are given $n$ general hyperplanes through the origin in $\R^{n-r}$,
and we replace each of them by a random parallel translate. The
resulting arrangement of $n$ affine hyperplanes in $\R^{n-r}$ creates
precisely $\binom{n-1}{r-1}$ bounded regions, as promised by 
the conjunction of  (\ref{eq:moebius2}) and (\ref{eq:moebius4}).
\hfill $\diamond$
\end{example}

The M\"obius number is important to us because it computes
the degree of the central curve of  the primal linear program (\ref{initial problem}).
See Theorem \ref{thm:main} below. 
We will take $r = d+1$ and $\mathcal{L}= \mathcal{L}_{A,{\bf c}}$ to be the linear space 
spanned by the rows of $A$ and the vector ${\bf c}$. 
The matroid $M(\L_{A,{\bf c}})$ we need there has rank 
$r=d+1$ and it is denoted  $M_{A,{\bf c}}$. 
We use the notation
\begin{equation}
\label{eq:moebius6}
|\mu(A,{\bf c})|\,\, := \,\, |\mu(M_{A,{\bf c}})| \,\, = \,\, |\mu(M(\mathcal{L}_{A,{\bf c}}))|. 
\end{equation}

Consider Example~\ref{transport23}, with $A$ the displayed
 $5 {\times} 6$-matrix of rank $d=4$, or the instance in
 Figure  \ref{fig:polygons}.  Here,
 $n=6$, $r=d+1=5$, and $M_{A,{\bf c}} = U_{5,6}$ is  the uniform matroid.
 Its M\"obius number equals  $|\mu(A,{\bf c})| = |\mu(U_{5,6})| = 5$.
 This number $5$ counts the bounded segments on
the vertical line on the left in  Figure  \ref{fig:polygons}. 
Note that the relevant matroid for  Example~\ref{transport23} is not, as one might expect,
 the graphic  matroid of  $K_{2,3}$. For higher-dimensional problems the  matroids 
$M_{A,{\bf c}}$ we encounter are typically non-uniform.

There are many other interpretations of the M\"obius invariant and the h-vector.
For example, a useful identity for computations is the h-vector as an evaluation of the Tutte polynomial 
$T_M(x,y)$ of the matroid $M$ (see \cite[Eq.~(7.15)]{bjornercomplex}  and the discussion  in \cite[\S 7.2]{bjornercomplex}): 
\begin{equation}
\label{eq:tutteH}
h_0 + h_1 z + h_2 z^2 + \cdots + h_{r-1} z^{r-1} \quad = \quad z^r \cdot T_M(1/z,0).
\end{equation}

\section{Equations defining the central curve}
\label{centralpathideal}

In this section we determine the prime ideal of the 
central curve of the primal linear program (\ref{initial problem}). 
As a consequence we obtain explicit formulas for the degree, 
arithmetic genus and Hilbert function of the projective closure of the primal 
central curve. These results resolve the problem stated
by Bayer and Lagarias at the end of  \cite[\S 11]{BL2}. 

Our ground field is $K$ will be  $\Q(A)({\bf b}, {\bf c})$. Here $\Q(A)$ denotes the
subfield of $\R$ generated by the entries of $A$ and $\Q(A)({\bf b}, {\bf c})$ is the 
rational function field generated by the coordinates $b_i$ and $c_j$ of the right 
hand side ${\bf b}$ and the cost vector ${\bf c}$. We assume that these coordinates 
are algebraically independent over $\Q(A)$. This is a formal way of ensuring that our 
algebraic results remain valid for generic values of $b_i$ and $c_j$ in $\R$.  
 
%This will ensure that all our algebraic results derived remain valid under almost 
% all other specializations $K \rightarrow \R$ of these coordinates to the field of real numbers.

Let $\L_{A,{\bf c}}$ be the subspace of $K^n$
spanned by the rows of $A$ and the vector ${\bf c}$.
We define the {\em central sheet} to be the coordinate-wise reciprocal
 $\L_{A,{\bf c}}^{-1}$ of that linear subspace.
In precise terms,  we define $\L_{A,{\bf c}}^{-1}$
to be the Zariski closure in the affine space $\mathbb{C}^n$ of the set
\begin{equation}\label{eq:centralsheet}
\biggl\{\, 
\left(
\frac{1}{u_1}, 
\frac{1}{u_2}, \ldots,
\frac{1}{u_n}
\right) \in \mathbb{C}^n \,\,:\,\,
(u_1,u_2,\ldots,u_n) \in \L_{A,{\bf c}}
\quad \text{and} \quad
u_i \not = 0 \,\,\, \text{for} \,\,\, i = 1,\ldots,n \biggr\}.\end{equation}

\begin{lemma} \label{lem:eins}
The Zariski closure of the primal central path
$\,\{ {\bf x}^*(\lambda) \,: \,\lambda \in \R_{\geq 0} \}\,$ is equal to the
intersection of the central sheet 
$\,\L_{A,{\bf c}}^{-1}\,$
with the affine-linear subspace defined by
$\,A {\bf x} = {\bf b} $.
\end{lemma}

\begin{proof}
We eliminate ${\bf s}, {\bf y}$ and $\lambda$ from the equations
$\, A^T {\bf y} - {\bf s}  = {\bf c}\,$ and  $\,x_i  s_i = \lambda \,$ as follows.
We first replace the coordinates of ${\bf s}$ by $s_i  = \lambda/x_i$. The linear system becomes
$\,A^T {\bf y} - \lambda {\bf x}^{-1} = {\bf c}$. This condition means that
${\bf x}^{-1} = (\frac{1}{x_1},\ldots,\frac{1}{x_n})^T$ lies in the
linear space $\L_{A,{\bf c}}$ spanned by ${\bf c}$ and the rows of $A$. The result
of the elimination says that ${\bf x}$ lies in the central sheet $\L_{A,{\bf c}}^{-1}$.
For ${\bf x}$ in the Zariski-dense set $\L_{A,{\bf c}}^{-1}\cap(\mathbb{C}^*)^n$, one can reconstruct 
values of $\lambda, {\bf y}, {\bf s}$ for which $({\bf x}, {\bf y}, {\bf s}, \lambda)$ is a solution to the equations
$A^T{\bf y}-{\bf s}={\bf c}$, $x_is_i =\lambda$. This shows that $\L_{A, {\bf c}}^{-1}$ is indeed the
 projection  of the set of these solutions
 onto the ${\bf x}$-coordinates.
\end{proof}

The linear space $\{A {\bf x} = {\bf b}\}$ has dimension
$n-d$, and we write $I_{A,{\bf b}}$ for its linear ideal.
The central sheet $\,\L_{A,{\bf c}}^{-1} \,$
is an irreducible variety of dimension $d+1$,
and we write $J_{A,{\bf c}}$ for its prime ideal.
Both $I_{A,{\bf b}}$ and $J_{A,{\bf c}}$ are ideals in $K[x_1,\ldots,x_n]$.
We argue the following is true:

\begin{lemma} \label{lem:zwei}
The prime ideal of polynomials that vanish on the central curve $\mathcal{C}$
is $\,I_{A,{\bf b}} + J_{A,{\bf c}}$.
The degree of both $\mathcal{C}$ 
and the central sheet $\,\L_{A,{\bf c}}^{-1}$ 
coincides with the M\"obius number $|\mu(A,{\bf c})|$.
\end{lemma}

\begin{proof} 
The intersection of the affine space $\{A {\bf x} = {\bf b}\}$
with the central sheet is the variety of
the ideal  $\,I_{A,{\bf b}} + J_{A,{\bf c}}$. This ideal is prime because ${\bf b}$ and ${\bf c}$
are generic over $\Q(A)$. The intersection is the central curve. In Proposition~\ref{prop:PS} we 
show that the degree of the central sheet is $|\mu(A,{\bf c})|$,
so here it only remains to show that this is the degree of the central curve as well. 
For a generic vector $({\bf b}, c_0) \in \R^{d+1}$, we consider the hyperplane 
arrangement induced by $\{x_i=0\}$ in the affine
 space $\{\binom{A}{{\bf c}}{\bf x} = \binom{{\bf b}}{c_0} \}$. 
The number of bounded regions of this hyperplane arrangement equals the
M\"obius number $|\mu(A,{\bf c})|$, as seen in (\ref{eq:moebius4}).
% Note that $|\mu(A,{\bf c})|$ does not depend on ${\bf c}$, since this vector is generic over~$\Q(A)$.

Each of these bounded regions contains a unique point maximizing $\sum_i \log|x_i|$. 
This point is the \textit{analytic center} of that region. Each such 
analytic center lies in $\L_{A,{\bf c}}^{-1}$, and thus on the central curve by Lemma~\ref{lem:eins}.
This shows that the intersection of the central curve with the plane $\{{\bf c}^T{\bf x} = c_0\}$ contains at least $|\mu(A,{\bf c})|$ points. 

B\'ezout's Theorem implies that the degree of a variety $V \subset \mathbb{C}^n$ is an upper bound for  the degree of its intersection $V\cap H$ with an affine subspace $H$, provided that $n+\dim(V\cap H) = \dim(V)+\dim(H)$. We use this theorem for two inequalities; first, that the degree of $\L_{A,{\bf c}}^{-1}$ bounds the degree of the central curve $\mathcal{C}$, and, second, that the degree of $\mathcal{C}$ bounds the number of its intersection points with $\{{\bf c}^T{\bf x} = c_0\}$. To summarize, we have shown:
\[
|\mu(A,{\bf c})|\;\; \leq\;\; \#(  \mathcal{C} \cap \{{\bf c}^T {\bf x}=c_0\} )\;\; \leq \;\; \deg(\mathcal{C})\;\; \leq  \;\;\deg( \L_{A,{\bf c}}^{-1}) \;=\; | \mu(A,{\bf c})|.\]
From this we conclude that $|\mu(A,{\bf c})|$ is the degree of the primal central curve $\mathcal{C}$.\end{proof}

At this point we are left with the problem of computing
the degree of the homogeneous ideal $J_{A,{\bf c}}$ and a set of generators. Luckily, this has already been done for us in the literature. The following proposition was proved by Proudfoot and Speyer \cite{PS}
and it refines an earlier result of Terao \cite{Ter}.  See also \cite{berget} for related results. 
The paper \cite{SU} 
begins the challenging task of extending these
results from linear programming to semidefinite programming.

\begin{proposition}[Proudfoot-Speyer \cite{PS}] 
\label{prop:PS}
The degree of the central sheet $\L_{A,{\bf c}}^{-1}$, regarded as a variety in complex projective space,
 coincides with the M\"obius number $|\mu(A,{\bf c})|$.
Its prime ideal $J_{A,{\bf c}}$ is generated by a universal Gr\"obner basis consisting of
all homogeneous polynomials
\begin{equation}
\label{eq:circuits}
 \sum_{i \in {\rm supp}(v)} v_i \cdot \prod_{j \in {\rm supp}(v)\backslash\{i\}} x_j , 
 \end{equation}
where $\sum v_i x_i$ runs over non-zero linear forms of minimal support
that vanish on $\L_{A,{\bf c}}$.
\end{proposition}

\begin{proof} 
The construction in \cite{PS} associates the ring
$K[x_1,\ldots,x_n]/J_{A,{\bf c}}$ to the linear subspace
$\L_{A,{\bf c}}$ of $K^n$. Theorem 4 of \cite{PS} says that the homogeneous polynomials
(\ref{eq:circuits}) form a universal Gr\"obner bases for $J_{A,{\bf c}}$. 
As argued in \cite[Lemma 2]{PS}, this means that the ring degenerates to the Stanley-Reisner ring of
the broken circuit complex  ${\rm Br}(M_{A,{\bf c}})$. Hence, by our discussion in Section 3,
or by \cite[Prop.~7]{PS},
the Hilbert series of $K[x_1,\ldots,x_n]/J_{A,{\bf c}}$ is the rational function (\ref{eq:FH}),
and the degree of $J_{A,{\bf c}}$ equals $|\mu(A,{\bf c})|$ as seen in (\ref{eq:moebius5}).
The ideal $J_{A,{\bf c}}$ is radical, since its initial ideal is square-free,
and hence it is prime because its variety $\mathcal{L}_{A,{\bf c}}^{-1}$ is irreducible.
\end{proof}

The polynomials in (\ref{eq:circuits}) correspond to the circuits
of the matroid $M_{A,{\bf c}}$, of which there are at most $\,\binom{n}{d+2}$.
%There is at most one circuit contained in each
%$(d+2)$-subset of $\{1,\ldots,n\}$, so their number is at most $\,\binom{n}{d+2}$. 
If the matrix $A$ is generic, then
$M_{A,{\bf c}}$ is the uniform matroid and every $(d+2)$-subset of $\{1,\ldots,n\}$ forms a circuit. 
In this case, by (\ref{eq:moebius2}),
 its M\"obius number equals 
$$ |\mu(A,{\bf c})| \,\,= \,\, \binom{n-1}{d}. $$
For arbitrary matrices $A$, this binomial coefficient furnishes an upper bound
on the M\"obius number $|\mu(A,{\bf c})|$.
We are now prepared to conclude with the main theorem of this section.
The analogous equations for 
the dual central curve are given in Proposition~\ref{thm:maindual} in Section~\ref{globalgeo}.

\begin{theorem} \label{thm:main} 
The degree of the primal  central curve of (\ref{initial problem}) is the M\"obius number
$|\mu(A,{\bf c})|$ and is hence at most $\binom{n-1}{d}$. The prime ideal of  polynomials that 
vanish on the primal central path is generated by the circuit polynomials
(\ref{eq:circuits}) and the $d$ linear polynomials in $A {\bf x} - {\bf b}$.
\end{theorem}

\begin{proof}
This is an immediate consequence of 
Lemmas \ref{lem:eins} and \ref{lem:zwei} and
Proposition \ref{prop:PS}. 
\end{proof}

It is convenient to write the circuit equations  (\ref{eq:circuits}) in the following
determinantal representation. Suppose that $A$ has format $d {\times} n$ and its
rows are linearly independent. Then the linear forms of minimal support that vanish 
on $\L_{A, {\bf c}}$ are the $(d+2)\times(d+2)$-minors of the $(d+2)\times n$ matrix 
${\tiny \begin{pmatrix} A \\ {\bf c} \\ {\bf x}\end{pmatrix} }$. 
This gives the following concise description of our prime ideal $J_{A,{\bf c}}$:
\vskip -0.3cm
\begin{equation}
\label{Jdet}
J_{A,{\bf c}} \quad = \quad
I_{{\rm num},d+2}
\begin{pmatrix} \,\, A \,\, \\ {\bf c} \\
\,\, {\bf x}^{-1} \,\,
\end{pmatrix}
\end{equation}
where ${\bf x}^{-1} = (x_1^{-1}, \ldots, x_n^{-1})$ and
the operator $I_{{\rm num},d+2}$ extracts 
the numerators of the $(d{+}2)\times (d{+}2)$-minors of the matrix.
For example, 
if the leftmost $(d+1)\times(d+2)$ submatrix of the matrix $\binom{A}{\bf c}$ has full rank $d+1$, 
then, as in \eqref{eq:sec2matrix}, one generator of the ideal $J_{A,{\bf c}}$ equals
\[
\; \det\begin{pmatrix} 
A_{1} \! & \! A_{2}  \! & \!  \hdots & \! A_{d+2} \\
c_{1} \! & \!  c_{2}  \! & \! \hdots \! & \! c_{d+2}\\
 x_{1}^{-1} & x_{2} ^{-1} & \hdots  & x_{d+2} ^{-1} \\
\end{pmatrix} \cdot \prod_{i\in \mathcal{I}} x_i,
\]
where $\mathcal{I}$ is the unique circuit contained in $\{1,2,\hdots, d+2\}$.
Note that there are $\binom{n}{d+2}$ such minors but they  need not be distinct
and some of them  may be zero.
%For example, one generator of $J_{A,{\bf c}}$ equals
%\[
%\; \det\begin{pmatrix} 
%A_{1} \! & \! A_{2}  \! & \!  \hdots & \! A_{d+2} \\
%c_{1} \! & \!  c_{2}  \! & \! \hdots \! & \! c_{d+2}\\
% x_{1}^{-1} & x_{2} ^{-1} & \hdots  & x_{d+2} ^{-1} \\
%\end{pmatrix} \cdot \prod_{i\in \mathcal{I}} x_i,
%\]
%where $\,\mathcal{I} $ is the lexicographically earliest circuit of the matroid $M_{A,{\bf c}}$. 
%Note that there are $\binom{n}{d+2}$ such minors but they need not be distinct.
%

\begin{example}
\label{transport23b}
Let $d=4,n=6$ and $A$ the matrix in Example \ref{transport23}. The linear ideal is
$$
I_{A,{\bf b}} \,\,\, = \,\,\, \langle  \, x_1 + x_2 + x_3 - b_1 \,, \, x_4 + x_5 + x_6 - b_2 \,,\, 
x_1 + x_4 - b_3 \,, \, x_2   + x_5 - b_4 \,\rangle. $$
The central sheet $\mathcal{L}_{A,{\bf c}}^{-1}$ is the quintic hypersurface 
whose defining polynomial is
\begin{equation} 
\label{sixbysixdet}
 f_{A,{\bf c}}({\bf x}) \quad = \quad
{\rm det}
\begin{small}
\begin{pmatrix}
1 & 1 & 1 & 0 & 0 & 0 \\
0 & 0 & 0 & 1 & 1 & 1 \\
1 & 0 & 0 & 1 & 0 & 0 \\
0 & 1 & 0 & 0 & 1 & 0 \\
c_1 & c_2 & c_3 & c_4 & c_5 & c_6 \\
x_1^{-1} & x_2^{-1} & x_3^{-1} & x_4^{-1} & x_5^{-1} & x_6^{-1} 
\end{pmatrix}
\end{small} 
\cdot x_1 x_2 x_3 x_4 x_5 x_6 .
\end{equation}
The primal central curve is the plane quintic defined by the ideal
$I_{A,{\bf b}} +  \langle f_{A,{\bf c}} \rangle $.
This ideal is prime for general choices of ${\bf b}$ and ${\bf c}$.
However, this may fail for special  values:
the quintic on the left in  Figure \ref{fig:hexagon}
 is irreducible but that on the right decomposes
into a quartic and a line.
For a concrete numerical example we set $ b_1 = b_2 = 3 $ and
$b_3 = b_4 = b_5 = 2$.  Then the transportation polygon $P$ is the
regular hexagon depicted in  Figure \ref{fig:hexagon}. Its vertices are
\begin{equation}
\label{hexagon}
\begin{small}
\begin{pmatrix}
0 & 1 & 2 \\
2 & 1 & 0 
\end{pmatrix},\,\,
\begin{pmatrix}
0 & 2 & 1 \\
2 & 0 & 1 
\end{pmatrix},\,\,
\begin{pmatrix}
1 & 0 & 2 \\
1 & 2 & 0 
\end{pmatrix},\,\,
\begin{pmatrix}
1 & 2 & 0 \\
1 & 0 & 2 
\end{pmatrix},\,\,
\begin{pmatrix}
2 & 0 & 1 \\
0 & 2 & 1 
\end{pmatrix},\,\,
\begin{pmatrix}
2 & 1 & 0 \\
0 & 1 & 2 
\end{pmatrix}
\end{small}.
\end{equation}
Consider the two transportation problems (\ref{initial problem}) given by 
$\,{\bf c} = $ \begin{small} $ \begin{pmatrix} 0 & 0 & 0 \\ 0 & 1 & 3 \end{pmatrix}\,$
\end{small}
 and
$\,{\bf c}' =$ \begin{small} $ \begin{pmatrix} 0 & 0 & 0 \\ 0 & 1 & 2 \end{pmatrix}$.
\end{small}
In both cases, the last matrix in (\ref{hexagon}) is the unique
optimal solution. Modulo the linear ideal $I_{A,{\bf b}}$ we
can write the quintics $f_{A,{\bf c}}$ and 
$f_{A,{\bf c}'}$ as polynomials in only two variables $x_1$ and $x_2$:
$$
f_{A,{\bf c}} \quad = \qquad
\begin{matrix}
3 x_1^4 x_2+5 x_1^3 x_2^2-2 x_1 x_2^4-3 x_1^4-22 x_1^3 x_2-15 x_1^2 x_2^2+8 x_1 x_2^3
+2 x_2^4 \\ +18 x_1^3+45 x_1^2 x_2-12 x_2^3-33 x_1^2-22 x_1 x_2+22 x_2^2+18 x_1-12 x_2,
\end{matrix}
$$
\vskip -0.2cm
$$
f_{A ,{\bf c}'} \quad = \qquad
\begin{matrix}
(x_2-1) \cdot (2 x_1^4 + 4 x_1^3 x_2 + x_1^2 x_2^2 - x_1 x_2^3 - 12 x_1^3 - 14 x_1^2 x_2 +
x_1 x_2^2  \\ \quad + x_2^3 + 22 x_1^2 + 10 x_1 x_2 - 5 x_2^2 - 12 x_1 + 6 x_2).
\end{matrix}
$$
Both quintics pass through all intersection points of the arrangement of six lines.
The cost matrix {\bf c} exemplifies the generic behavior, when
the quintic curve is irreducible. On the other hand,
the central path for ${\bf c}'$ is a segment on the horizontal line $x_2  = 1$
 in Figure \ref{fig:hexagon}. \hfill $\diamond$
\end{example} 

\begin{remark} When ${\bf b}$ or ${\bf c}$ is not generic, various aspects of the above analysis break down.  If ${\bf b}$ is not generic, then the hyperplane arrangement $\{x_i=0\}_{i\in [n]} \subset \{A{\bf x}={\bf b}\}$ may not be simple, that is, it may have a vertex at which more than $n-d$ hyperplanes meet. This vertex will maximize ${\bf c}^T{\bf x}$ over more than one adjoining region of the arrangement. In particular, the central curve passes through this vertex more than once and is singular at this point. 

If the cost function ${\bf c}$ is maximized at a (non-vertex) face of a region of the hyperplane arrangement $\{x_i=0\}_{i\in [n]} \subset \{A{\bf x}={\bf b}\}$, then the central curve meets this face in its analytic center and does not pass through any of the vertices of the hyperplane arrangement contained in the affine span of this face. For example, see Figure~\ref{fig:dtz6}. Another potential problem is that for non-generic ${\bf c}$ the curve defined by the equations of Theorem~\ref{thm:main} may be reducible, as happens for the cost vector ${\bf c'}$ in Example~\ref{transport23b}. The central curve will then be whatever component of these solutions passes through the region of interest. In particular, its degree and equations are no longer independent of the sign conditions on 
${\bf x}$.  Fortunately, the M\"obius number $|\mu(A,{\bf c})|$ is always an upper bound for the degree
of the central curve.
\hfill $\diamond$
 \end{remark}

In the remainder of this section we consider the question of what happens to the 
central sheet, and hence to the central path, when the cost function  ${\bf c}$ 
degenerates to one of the unit vectors $e_i$. Geometrically this means that the cost vector becomes normal to one of the constraint hyperplanes, and the curve reflects this by breaking into irreducible components. 
In addition to the nice geometry, our interest in these degenerations comes from 
the observation that they seem to produce curves with high curvature. 
What follows is independent of the rest of the paper and can be skipped upon first reading. 

To set up our degeneration in proper algebraic terms, we work over the field
$K \{\!\{ t \}\!\} $ of Puiseux series over the field $K = \Q(A)({\bf b},{\bf c})$ that was used above.
The field  $K \{\!\{ t \}\!\} $ comes with a natural $t$-adic valuation. Passing to the
special fiber represents the process of letting the parameter $t$ tend to $0$.
Our cost vector ${\bf c}$ has its coordinates in the Puiseux series field:
\begin{equation}
\label{CwithT}
 {\bf c}\,\, = \,\, \bigl(t^{w_1} , t^{w_2}, \ldots, t^{w_{n-1}} , 1 \bigr) 
 \end{equation}
Here $w_1 > w_2 > \cdots > w_{n-1} >  0 $ are any rational numbers.
We are interested in the special fiber of the central sheet
$\L_{A,{\bf c}}^{-1}$. This represents the limit of the central sheet
 as $t$ approaches $0$.
This induces a degeneration of the central curve
$\,\L_{A,{\bf c}}^{-1} \,\cap \,\{A {\bf x} = {\bf b} \}$.
We wish to see how, in that limit, the central curve breaks into irreducible curves in the
affine space $\{A {\bf x} = {\bf b} \}$.

The ideal defining the special fiber of $J_{A,{\bf c}}$ is denoted
$\,{\rm in}(J_{{A},{\bf c}}) = {J_{A,{\bf c}}}|_{t=0}$.
By a combinatorial argument as in \cite{PS}, 
 the maximal minors in (\ref{Jdet}) have the
Gr\"obner basis property for this degeneration.
Hence we obtain the prime ideal of the flat family by simply dividing each such minor by 
a non-negative power of $t$. 
This observation implies
the following result:

\begin{theorem} \label{thm:degen}
The central sheet $\L_{A,{\bf c}}^{-1}$ degenerates into a reduced union
of central sheets of smaller linear programming instances. More precisely, the ideal
$\,{\rm in}(J_{{A},{\bf c}}) \,$ is radical, and it has the following 
representation as an intersection of ideals that are prime when $A$ is generic:
\begin{equation}
\label{eq:pieces}
 {\rm in}(J_{{A},{\bf c}}) \quad = \quad
\bigcap_{i=d}^{n-1} \biggl( I_{{\rm num},d+1} \begin{pmatrix} A_1 & A_2 & \cdots & A_i \\
x_1^{-1} \! &\! x_2^{-1} \! & \cdots & \! x_i^{-1} 
\end{pmatrix}
\,+\, \langle x_{i+2}, x_{i+3}, \ldots , x_n \rangle \biggr)
\end{equation}
\end{theorem}

\begin{proof}[Proof sketch] 
The Gr\"obner basis property says that ${\rm in}(J_{A,{\bf c}})$
is generated by the polynomials obtained from the maximal minors
of (\ref{Jdet}) by dividing by powers of $t$ and then setting $t$ to zero.
The resulting polynomials factor, and this factorization shows that
they lie in each of the ideals on the right hand side of (\ref{eq:pieces}).
Conversely, each element in the product of the ideals on  the right hand side 
is seen to lie in
 ${\rm in}(J_{A,{\bf c}})$. To complete the proof, it then suffices to
 note that ${\rm in}(J_{A,{\bf c}})$ is radical because its
 generators form a square-free Gr\"obner basis.
 \end{proof}

\begin{example}
Let $n=6$ and $d=3$. The matrix $A$ 
might represent the three-dimensional
Klee-Minty cube. The decomposition of the initial ideal
in (\ref{eq:pieces}) has three components:
$$ 
{\rm in}(J_{A,{\bf c}}) \,\,= \,\,
\langle x_5 , \,x_6 \rangle \,\,\cap \,\,
\langle {\rm det} \!
\begin{pmatrix} x_1 A_1 \! & \! x_2 A_2 \! & \!  x_3 A_3 \! & \! x_4 A_4 \\
                             1 & 1 & 1 & 1 \\
\end{pmatrix} \! , x_6  \,\rangle  \,\,\cap \,\,
I_{{\rm num},4} \!\begin{pmatrix}
A_1 & A_2 & A_3 & A_4 & A_5 \\
\! x_1^{-1} \! & \! x_2^{-1} \! & \! x_3^{-1} \! & \! x_4^{-1} \! & \! x_5^{-1} \! 
\end{pmatrix}\! .  $$                                                          
For general $A$, the ideal $J_{A,{\bf c}}$ defines an irreducible curve of 
degree $10$, namely the central path, in each of the $3$-planes $\{A{\bf x} = {\bf b}\}$. The
three curves in its degeneration above are irreducible of degrees
$1$, $ 3$ and $6$ respectively. The first  is one of the
lines in the arrangement of six facet planes, the second curve
is the central path inside the facet defined by $x_6 = 0$,
and the last curve is the central path of the polytope
obtained by removing that facet. \hfill $\diamond$
\end{example}

In general, we can visualize the degenerated central path in the
following geometric fashion. We first flow from the analytic center 
of the polytope to the analytic center of its last facet. Then we iterate 
and flow from the analytic center of the facet to the analytic center of its last facet,
which is a ridge of the original polytope. Then we continue inside  that ridge, etc.

\section{The Gauss Curve of the Central Path}
\label{gaussc}

The total curvature of the central path is an important quantity for the estimation of 
the running time of interior point methods in linear programming 
\cite{DMS, MonteiroTsuchiya, Sonnevendetal,VavasisYe,ZhaoStoer}.
In this section we relate the algebraic framework developed so far to 
the problem of bounding the total curvature. The relevant geometry 
was pioneered by Dedieu, Malajovich and Shub \cite{DMS}. Following their
approach, we consider the {\em Gauss curve} associated with the primal central path.
The Gauss curve is the image of the central curve under the Gauss map, and its
 arc length is precisely the total curvature
of the central path. Moreover, the  arc length of the Gauss curve can be bounded in terms
of its degree. An estimate of that degree, via the multihomogeneous B\'ezout Theorem, was
the workhorse in \cite{DMS}. Our main result here is a more precise bound, 
in terms of matroid invariants, for the degree of the Gauss curve of the primal central curve.
As a corollary we obtain a new upper bound on its total curvature.

 We begin  our investigation by reviewing definitions  from elementary differential geometry.
Consider an arbitrary curve $\,[a,b] \rightarrow \R^m, \,t \mapsto f(t),\,$
whose parameterization is twice differentiable and whose
derivative $f'(t)$ is a non-zero vector for all parameter values $t \in [a,b]$.
This curve has an associated {\em Gauss map} into the unit sphere $S^{m-1}$, which is defined as
$$\gamma \,:\,[a,b] \rightarrow S^{m-1} \,,\,\,t \mapsto \frac{f'(t)} {||f' ( t )||}.$$
The image $\gamma = \gamma([a,b])$ of the Gauss map in $S^{m-1}$ is called the
{\em Gauss curve} of the given curve $f$. In our situation, the curve $f$
is algebraic, with known defining polynomial equations, and it makes sense
to consider the {\em projective Gauss curve} in complex projective space $\PP^{m-1}$.
By this we mean the Zariski closure of the image of the Gauss curve under
the double-cover map $S^{m-1} \rightarrow \PP^{m-1}$. If $m=2$, so that
$\C$ is a non-linear plane curve, then the Gauss curve traces out several
arcs on the unit curve $S^1$, and the projective Gauss curve is 
the entire projective line $\PP^1$. Here, the
line $\PP^1$ comes with a natural multiplicity, to be derived in
Example \ref{ex:sec5intheplane}.

If $m=3$ then the Gauss curve lies on the unit sphere $S^2$ and the
projective Gauss curve lives in the projective plane $\PP^2$. Since a curve in
$3$-space typically has parallel tangent lines, the Gauss curve is here
expected to have singularities, even if $f$ is a smooth curve.

The \emph{total curvature} $K$ of  our curve $f$ is defined to be the arc length of its
associated Gauss curve $\gamma$; see~\cite[\S 3]{DMS}. This quantity admits the following
expression as an integral:
\begin{equation}\label{curvaturedef} K \,\, := \, \int_{a}^{b} || \frac{d\gamma(t)}{dt} || dt.\end{equation} 

The {\em degree} of the Gauss curve $\gamma(t)$ is defined as the maximum number 
of intersection points, counting multiplicities, with any hyperplane in $\R^m$,
or equivalently, with any equator in $S^{m-1}$.
This (geometric) degree is bounded above by the (algebraic) degree of
the projective Gauss curve in $\PP^{m-1}$. The latter can be
computed exactly, from any polynomial representation of $\C$,
using standard methods of computer algebra.  Throughout this section,
by {\em degree} we mean the degree of the image of $\gamma$ in $\PP^{m-1}$
multiplied by the degree of the map that takes $\C$ onto $\gamma(\C)$.
From now on we use the notation ${\rm deg}(\gamma(\C))$ for that number.

\begin{proposition}
\label{prop:curva}
{\rm\cite[Corollary 4.3]{DMS}}
The total curvature of any real algebraic curve $\C$ in $\R^m$ is bounded above by
$\pi $ times the degree of its projective Gauss curve in $\PP^{m-1}$. In symbols,
$$ K \,\,\leq \,\, \pi \cdot {\rm deg}(\gamma(\C)) .$$
\end{proposition}

\begin{remark}
In higher dimensions, the degree of the map from $\C$ onto $\gamma(\C)$ is typically equal to $1$,
in which case our definition of ${\rm deg}(\gamma(\C))$ is exactly that
used in \cite{DMS}. However, that the extra factor is needed 
can be seen by considering the case $m = 2$ of non-linear plane curves:
the Gauss curve $\gamma(\C)$ is just $\PP^1$, but it has
 a non-reduced structure coming from the map.
 \hfill $\diamond$
\end{remark}

We now present our main result in this section, which concerns 
the degree of the projective Gauss curve $\gamma(\C)$,
when $\C$ is the central curve of a linear program in primal formulation.
As before, $A$ is an arbitrary real matrix of rank $d$
 having $n$ columns, but the cost vector ${\bf c}$ and the right hand side ${\bf b}$ are generic
 over $\Q(A)$.
 The curve $\C$ lives in an $(n-d)$-dimensional affine subspace of $\R^n$,
which we identify with  $\R^{n-d}$, so that $\gamma(\C)$ is a curve in $\PP^{n-d-1}$.

Let $M_{A,{\bf c}}$ denote the matroid of rank $d+1$ on the ground set
$[n]$ associated with the matrix $\binom{A}{{\bf c}}$.
We  write $(h_0,h_1,....,h_d)$ for the h-vector of the broken circuit complex of $M_{A,{\bf c}}$,
as defined in (\ref{eq:FH}). In the generic case, $M_{A,{\bf c}} = U_{d+1,n}$ is the uniform
matroid, as in Example~\ref{ex:uniform}.
%the maximal simplices in ${\rm Br}(M_{A,{\bf c}})$ are
%$\{1,j_1,\ldots,j_{d}\}$ where $2\leq j_1 < \cdots < j_{d} \leq n$.
In this case, the coordinates of the h-vector are $\,h_i =  { n-d+i-2 \choose i }$.
For special matrices $A$, this simplicial complex gets
replaced by a pure shellable subcomplex of the same dimension,
so the h-vector (weakly) decreases in each entry. Hence, the following always holds:
\begin{equation}
\label{eq:h_ibound}
 h_i \,\, \leq \,\,  { n-d+i-2 \choose i } \qquad {\rm for} \,\,\, i = 0,1,\ldots,d. 
 \end{equation}
 As indicated,  this inequality holds with  equality when $M_{A,{\bf c}}$ is the uniform matroid.

\begin{theorem} \label{thm:gaussmain}
The degree of the projective Gauss curve of
the primal central curve $\C$ satisfies
\begin{equation}
\label{eq:gausstruth}  {\rm deg}(\gamma(\C)) \,\, \leq \,\,\,  2 \cdot \sum_{i=1}^d i \cdot h_i . \qquad \quad
\end{equation}
In particular, we have the following upper bound which is tight for generic matrices $A$:
\begin{equation} 
\label{eq:h_i}
 \qquad  {\rm deg}(\gamma(\C)) \quad \leq \quad  2 \cdot (n-d-1)\cdot \binom{n-1}{d-1}.
 \end{equation}
\end{theorem}

The difference between the bound in (\ref{eq:gausstruth}) and the
degree of $\gamma(\mathcal{C})$ can be explained in terms of
singularities the curve $\mathcal{C}$ may have on the hyperplane
at infinity. The relevant algebraic geometry will be seen in
  the proof of Theorem \ref{thm:gaussmain},
which we shall present after an example.

\begin{example} \label{ex:gauss1}
In the following two instances we have $d=3$ and $n=6$. 
\begin{enumerate}

\item First assume that $A$ is a generic $3 \times 6$-matrix. The arrangement of six
facet planes creates $10$ bounded regions.  The primal
central curve $\C$ has degree $\binom{6-1}{3} = 10$. It passes through the
 $\binom{6}{3}=20$ vertices of the arrangements. In-between it
 visits the $10$ analytic centers of the bounded regions. Here the curve $\mathcal{C}$
 is smooth and its genus is $11 $. This number is seen from
  the formula (\ref{eq:genusformula}) below.
The corresponding Gauss curve in $\mathbb{P}^2$ has degree
$2\cdot 10 + 2 \cdot {\rm genus}(\C) - 2 = 40$, as given by
 the right hand side of (\ref{eq:h_i}). 
Hence the total curvature of the central curve $\C$ is bounded above by $\,40 \pi$. 

\item Next consider the Klee-Minty cube in $3$-space. Normally, it is given by the constraints
$$
0 \leq z_1 \leq 1 \,, \,\,\,
\epsilon z_1 \leq z_2 \leq 1 - \epsilon z_1, \quad \hbox{and} \quad
\epsilon z_2 \leq z_3  \leq 1 - \epsilon z_2.
$$
To see this in a primal formulation (\ref{initial problem}), we use $z_1,z_2,z_3$ to parametrize the affine space $\{A{\bf x}={\bf b}\}$. 
The facets of the cube then correspond to the intersection of the coordinate hyperplanes with this affine space. 
This is given by the matrices
$$ \binom{A}{{\bf c}}  \quad = \quad
\begin{pmatrix}
  1 & 1 & 0 & 0 & 0 & 0 \\
   2 \epsilon & 0 & 1 & 1 & 0 & 0 \\
   2 \epsilon^2  & 0 & 2 \epsilon & 0 & 1 & 1 \\
   c_1 & c_2 & c_3 & c_4 & c_5 & c_6 
   \end{pmatrix}  \;\;\; \text{ and }\;\;\; {\bf b} = \begin{pmatrix} 1\\1\\1\end{pmatrix}  .$$
Here $\epsilon$ is a small positive real constant.
The above $4 {\times} 6$-matrix is not generic, and its associated matroid
$M_{A,{\bf c}}$ is not uniform. It has exactly one non-basis, and so the h-vector equals
$\,(h_0,h_1,h_2,h_3) = (1,2,3,3) $.
The central curve $\C$ has degree $\,\sum_{i=0}^3 h_i = 9$. In the
coordinates used above, the curve is defined by
the $5 {\times} 5$-minors of the $5 {\times} 6$-matrix
which is obtained from the $4 {\times} 6$-matrix $\binom{A}{{\bf c}}$ by
adding one row consisting of reciprocal facet equations:
$$ \qquad \qquad \bigl(\, z_1^{-1},\, (1-z_1)^{-1},\, (z_2- \epsilon z_1)^{-1}, \,
(1-z_2- \epsilon z_1)^{-1}, \,(z_3- \epsilon z_2)^{-1}, \,(1-z_3- \epsilon z_2)^{-1} \, \bigr).$$
According to Theorem \ref{thm:gaussmain}, the
degree of the Gauss curve $\,\gamma(\C)\,$ in $\,\PP^2\,$ is bounded above
by $\,2\sum_{i=1}^3 i \cdot h_i = 34 $.  A computation using {\tt Macaulay2} \cite{M2} reveals that 
${\rm degree}(\gamma(\C)) = 32$ and that the total curvature is bounded by $32\pi$. \hfill $\diamond$
\end{enumerate}
\end{example}

\medskip

\begin{proof}[Proof of Theorem \ref{thm:gaussmain}]
For the proof we shall use the {\em generalized Pl\"ucker formula for curves}:
\begin{equation}
\label{eq:plucker}
 {\rm deg}(\gamma(\C)) \,\, = \,\, \, 2 \cdot {\rm deg}(\C) + 2 \cdot {\rm genus}(\C) - 2 - \kappa.
 \end{equation}
  The formula in (\ref{eq:plucker}) is obtained from \cite[Thm.~(3.2)]{piene}
 by setting $m=1$ or from \cite[Eq.~(4.26)]{griff} by setting $k=0$.
 The quantity $\kappa$ is a non-negative integer and it measures the singularities of the curve $\C$.
We have $\kappa = 0$ whenever the projective curve $\C$ is smooth, and this happens in our application
when $M_{A,{\bf c}}$ is the uniform matroid. In general, we may have singularities
at infinity because here the real affine curve $\C$ has to be replaced by its closure
 in complex projective space $\PP^{n-d} $, which is the projectivization of the
 affine space defined by $A {\bf x} = {\bf b}$. The degree and genus on the
 right hand side of (\ref{eq:plucker}) refer to that projective curve  in $\PP^{n-d}$. 
 
The references above actually give the degree 
of the {\em tangent developable} of the projective curve $\C$, but we see that this equals the degree of the Gauss curve. 
The tangent developable is the surface obtained by taking the union of all tangent lines at points in $\C$.
The projective Gauss curve $\gamma(\C)$ is obtained from the
tangent developable by intersecting it with a hyperplane,
namely, the hyperplane at infinity, representing the directions of lines.

In the formula (\ref{eq:plucker}), the symbol ${\rm genus}(\C )$ refers to the
arithmetic genus of the curve. We shall now compute this arithmetic genus for
primal central curve $\C$.
For this we use the formula for the 
Hilbert series of the central sheet due to Terao, in Theorem 1.2 on page 551 of \cite{Ter}.
See the recent work of Berget \cite{berget} for
a nice proof of a more general statement.

As seen in the proof of Proposition \ref{prop:PS},
the Hilbert series of the coordinate ring of the central sheet equals
$$ \frac{h_0 + h_1 z + h_2 z^2 + \cdots + h_d z^d}{(1-z)^{d+2}}. $$
The central curve $\C$ is obtained from the central sheet by intersection with a general
linear subspace of dimension $n-d$. 
The (projective closure of the) central sheet is arithmetically Cohen-Macaulay
since it has a flat degeneration to a shellable simplicial complex, as shown by
 Proudfoot and Speyer \cite{PS}. We  conclude that
the Hilbert series of the central curve $\C$ is 
$$
\frac{h_0 + h_1 z + h_2 z^2 + \cdots + h_d z^d}{(1-z)^2} \,=\,
\sum_{m \geq d} \biggl[ (\sum_{i=0}^d h_i) \cdot m \,+ \sum_{j=0}^d (1-j) h_j \biggr] z^m \,\,+ \,\, O(z^{d-1}).
$$
The parenthesized expression is the Hilbert polynomial of the projective curve $\C$.
The degree of $\C$ is the coefficient of $m$, and using \eqref{eq:moebius5} we again find this to be the M\"obius number:
%relating the degree and  M\"obius number:
$$ {\rm degree}(\C) \,\,=\,\, |\mu(A,{\bf c})| \,\, = \,\, \sum_{i=0}^d h_i .$$
The arithmetic genus of the curve $\C$ is one minus the constant term of its Hilbert polynomial:
\begin{equation}
\label{eq:genusformula} {\rm genus}(\C) \,\, = \,\, 
1 - \sum_{j=0}^d (1-j) h_j .
\end{equation}
We now see that our assertion (\ref{eq:gausstruth}) follows  directly from 
the generalized Pl\"ucker formula~(\ref{eq:plucker}).

For fixed $d$ and $n$, the degree and genus of $\C$ are maximal 
when the matrix $A$ is generic. In this case, $h_i$ equals the right hand side of
(\ref{eq:h_ibound}), and we need to sum these binomial coefficients times two.
Hence, our second assertion (\ref{eq:h_i}) follows from the identify
$$ \sum_{i=0}^d \,i \cdot \binom{n-d+i-2}{i} \quad = \quad 
 (n-d-1)\cdot \binom{n-1}{d-1}. $$
 This completes the proof of Theorem \ref{thm:gaussmain}. \end{proof}

\begin{example} \label{ex:sec5intheplane}
Let $d=n-2$ and suppose $A$ is generic. The primal central curve $\C$ is a plane curve.
Our h-vector equals $(1,1, ...,1)$. The proof of Theorem \ref{thm:gaussmain} reveals  that
the degree of $\C$ is $d+1=n-1$ and the genus of $\C$ is ${d \choose 2}$.
The Gauss curve $\gamma(\C)$ is the projective line $\PP^1$, but
regarded with multiplicity $\,{\rm deg}(\gamma(\C)) = (d+1)d$. This number is
the degree of the projectively dual curve $\C^\vee$. The identity
 (\ref{eq:plucker}) specializes to the  Pl\"ucker formula for plane curves, 
 which expresses the degree of $\C^\vee $ in terms of the  degree and the
  singularities of $\C$. \hfill $\diamond$
 \end{example}
 
We close this section by showing
how to compute the Gauss curve for a non-planar instance.

 \begin{example}
 Let $n=5, d=2$ and $A  = \begin{pmatrix} 1 & 1 & 1 & 1 & 1 \\ 0 & 1 & 2 & 3 & 4 \end{pmatrix}$.
 The primal central curve has degree $6$ and its equations
 are obtained by clearing denominators in the $4 {\times} 4$-minors of
 $$ \begin{pmatrix} 1 & 1 & 1 & 1 & 1 \\ 0 & 1 & 2 & 3 & 4 \\
c_1 & c_2 & c_3 & c_4 & c_5 \\
(z_1 - g_1)^{-1} &  (-2 z_1 + z_2 - g_2)^{-1} &  (z_1 - 2 z_2 + z_3 - g_3)^{-1} &
     (z_2 - 2 z_3 - g_4)^{-1} &  (z_3 - g_5)^{-1}
     \end{pmatrix}. $$
The $c_i$ and $g_j$ are random constants representing the
cost function and right hand side of (\ref{initial problem}).
To be precise, the vector ${\bf g} = (g_1,g_2,g_3,g_4,g_5)^T$
satisfies $A  {\bf g} = {\bf b}$ as in  Sections \ref{dual_and_average} and \ref{globalgeo}.
The linear forms in ${\bf z}=(z_1,z_2,z_3)$ come from the change of coordinates $B^{T}{\bf z}-{\bf x}={\bf g}$ 
where $B$ is a $3\times 5$ matrix whose rows span the kernel of $A$.
Writing $\,{\tt I}\,$ for the ideal of these polynomials, the following one-line command
in the computer algebra system {\tt Macaulay2} \cite{M2} computes
the defining polynomial of the Gauss curve in $\PP^2$:
\begin{footnotesize}
\begin{verbatim}
  eliminate({z1,z2,z3},I+minors(1,matrix{{u,v,w}}*diff(matrix{{z1},{z2},{z3}},gens I)))
\end{verbatim}
\end{footnotesize}
The output is a homogeneous polynomial of degree $16$ in the
coordinates $u,v,w$ on $\PP^2$.
Note that ${\rm deg}(\gamma(\C)) = 16$ is consistent with Theorem \ref{thm:gaussmain} because
$h = (h_0,h_1,h_2) = (1,2,3)$. \hfill $\diamond$
\end{example}

\section{The Primal-Dual Formulation and the Average Total Curvature}
\label{dual_and_average}

This short section offers a dictionary that translates between the
primal and the dual central curve.   We begin by offering an algebraic representation
of the primal-dual central curve that is more symmetric
than that given in the Introduction. This will allow us to take our main results so far,
previously only stated for the primal LP, and show that they hold verbatim also for the dual LP. 
As primary application, we then derive a matroid-theoretic refinement of the
Main Theorem in \cite{DMS} on the average total curvature of the central path.

Let $\L_{A}$ denote the row space of the matrix $A$
and $\L_{A}^\perp$ its orthogonal complement in $\R^n$.
We also fix a vector ${\bf g} \in \R^n$ such that
$A{\bf g} = {\bf b}$. By eliminating 
${\bf y}$ from the system (\ref{eq:central1}) in
Theorem \ref{introthm}, we see that the primal-dual
central path $({\bf x^*}(\lambda),{\bf s^*}(\lambda))$ has the following 
symmetric description:
\begin{equation}
\label{eq:primaldual1}
   {\bf x} \in \L_{A}^\perp + {\bf g}\;,\,\,\,\;\;  {\bf s} \in \L_{A} + {\bf c}\,\,\,\;\;    \,\text{and} \,\,\, \;\; x_1 s_1 = x_2 s_2 = \dots = x_n s_n
   = \lambda .
   \end{equation}
The implicit (\textit{i.e.}~$\lambda$-free) representation of the primal-dual central curve
is simply obtained by erasing the very last equality ``$= \lambda$''
in (\ref{eq:primaldual1}). Its prime ideal is generated by the quadrics
$x_i s_i - x_j s_j$ and the affine-linear equations defining $\L_{A}^\perp + {\bf g}$ in ${\bf x}$-space and $\L_{A}+ {\bf c}$ 
in ${\bf s}$-space.

The symmetric description of the central path  in (\ref{eq:primaldual1}) lets us 
write down the statements from Section \ref{centralpathideal} for the dual version.
 For example, we derive equations for the dual central curve in ${\bf s}$-space or ${\bf y}$-space 
 as follows. Let $B$ be any $(n-d) \times n$ matrix whose rows span the kernel of $A$.
 In symbols, $\L_B=\L_A^{\perp}$.
%  Let ${\bf g}\in \R^n$ be any vector satisfying $A {\bf g} = {\bf b}$. 
 The (dual) central curve in ${\bf s}$-space is obtained by intersecting the $d$-dimensional affine space
$\, \L_{A}+ {\bf c}\, = \{\,{\bf s} \in \R^n \,:\,B{\bf s} = B {\bf c}\,\}\,$ 
with the  central sheet $\mathcal{L}_{B,{\bf g}}^{-1}$ in (\ref{eq:centralsheet}).  
To obtain the central curve in ${\bf y}$-space, we 
substitute $s_i = \sum_{j=1}^d a_{ji} y_j - c_i$ in the equations defining
$\mathcal{L}_{B,{\bf g}}^{-1}$.
This gives dual formulations of Theorems \ref{thm:main} and \ref{thm:gaussmain}: 

\begin{corollary} \label{thm:maindual}
The degree of the dual central curve of (\ref{dual problem slack}) equals the M\"obius number
$|\mu(B,{\bf g})|$ and is hence at most $\binom{n-1}{d-1}$. The prime ideal of  polynomials that
vanish on the central path is generated by the circuit polynomials (\ref{eq:circuits}), but now associated 
with the space generated by the rows of $B$ and the vector ${\bf g}$, and the $n-d$
linear equations in ${\bf s}$ given by $B {\bf s} = B {\bf c}$. \end{corollary}

\begin{corollary} \label{thm:gaussmaindual}
The degree $\deg(\gamma(\C))$ of the Gauss image of the dual central curve $\C$ 
is at most $2 \cdot \sum_i i \cdot h_i$, where $h=h({\rm Br}(M_{B,{\bf g}}))$.
This implies the bound $\,\deg(\gamma(\C)) \leq 2\cdot(d-1)\cdot\binom{n-1}{d}$.
\end{corollary}

\begin{remark}
Theorem~\ref{thm:gaussmain} and Corollary \ref{thm:gaussmaindual} give a strengthening of Theorem 1.1 of
\cite{DMS}. Megiddo and Shub \cite{megiddoshub} proved the lower bound  $ (d-1)\pi/2$
for  the total curvature of the central path of a $d$ variable LP with $d+1$ constraints.
For such instances, our upper bound is tight up to a constant. Our main contribution is that we
adjust the upper bound to the specific  matroid of the constraint matrix $A$.  \hfill $\diamond$
\end{remark}

Dedieu, Malajovich and Shub \cite{DMS} investigated the average curvature of the 
central paths over all the bounded regions of the associated hyperplane arrangement for the linear program 
and show this to be at most $2\pi d$ when $A$ is generic. Combining \eqref{eq:moebius4}, 
Proposition~\ref{prop:curva},  Theorem~\ref{thm:gaussmain}, and their dual versions seen in this section, we 
obtain the following statement:

 \begin{theorem}\label{avgCurv} 
    The average total curvature of the primal central path over the bounded regions of the
  hyperplane arrangement $\{x_i=0\}_{i\in [n]}$ in
  the affine space $\,\L_{A}^\perp + {\bf g}\,$
     is at most 
  \begin{equation}
  \label{eq:2piPrimal}
   2 \pi \cdot \frac{ \sum_{i=0}^d i \cdot h_i(A,{\bf c})}{|\mu(A)|} .
   \end{equation}
  If $A$ is a generic $d\times n$ matrix, then this average total curvature is at most $\;2 \pi (n-d-1)  $.

 Similarly, the  average total curvature of the dual central path over the bounded regions of the
  hyperplane arrangement $\{s_i=0\}_{i\in [n]}$ in
  the affine space $  \,   \L_{A} + {\bf c} \,$ is at most
 \begin{equation}
  \label{eq:2piDual}
  2 \pi \cdot \frac{ \sum_{i=0}^d i \cdot h_i(B,{\bf g})}{|\mu(B)|} .
  \end{equation}
  In the generic case, this average total curvature is bounded above by $ \,2 \pi (d-1)$.
 \end{theorem}
 
 Note that the h-vector entries in the numerator of (\ref{eq:2piPrimal})
refer to the rank $d+1$ matroid $M_{A,{\bf c}}$ while the
M\"obius number in the denominator refers to the rank $d$ matroid $M_A$.
 Likewise, the  h-vector entries in the numerator of (\ref{eq:2piDual})
 refer to the rank $n-d+1$ matroid $M_{B,{\bf g}}$ while the
M\"obius number in the denominator refers to the rank $n-d$ matroid $M_B$.
 
 \begin{example} \label{ex:backtothesnake} 
We consider two examples to demonstrate the finer bounds provided by Theorem \ref{avgCurv}.
First, recall the primal LP formulation of the Klee-Minty Cube in Example \ref{ex:gauss1}.
We had already seen that our formulas predict the total curvature is bounded by $32\pi$.
We can calculate $|\mu(A)|$ for the $3 \times 6$ defining matrix. The matroid in question is
not uniform (it only has 14 bases). We recover $\mu(A)$ through the computation and evaluation
of the Tutte polynomial as shown in \ref{eq:tutteH}.  In this case we obtain $|\mu(A)|=5$ 
thus Formula~\ref{eq:2piPrimal} says the average curvature is no more than $\frac{32 \pi}{5}$. This is a better
bound than the generic case formula.

Second, recall the DTZ snake of Example \ref{ex:snake} and Figure~\ref{fig:dtz6} 
which is given in dual LP form. The curve has degree four and its projective closure
$\mathcal{C}$ is smooth in $\PP^2$. So, we have ${\rm deg}(\gamma(\C)) = 12$, and Proposition~\ref{prop:curva} gives an 
upper bound $12 \pi$ on the total curvature of the entire central curve in $\R^2$. The M\"obius number $|\mu(B)|$ is in this case
10, thus we get from formula~(\ref{eq:2piDual}) that the average curvature is $1.2 \pi$. Once more, our bound is better than the
generic case prediction. It is interesting to see the large discrepancy between the curvature within 
a single feasible region (calculated in Example \ref{ex:snake}) and the average total curvature. \hfill $\diamond$
\end{example}
 
\section{Global Geometry of the Central Curve}
\label{globalgeo}

In this section we return to the central path in its 
primal-dual formulation, and we study several geometric properties.
We shall show how the central curve connects the
vertices of the hyperplane arrangement with the
analytic centers of its bounded regions. This picture
behaves well under duality, as the vertices of the two arrangements are in 
natural bijection.

For an algebraic geometer, it is natural to 
replace each of the affine spaces in (\ref{eq:primaldual1}) by a complex
projective space  and to study the closure $\mathcal{C}$ of
the central curve in  $\PP^n \times \PP^n$.
Algebraically, we use homogeneous coordinates
$ [x_0:x_1:\cdots:x_n]$ and $[s_0:s_1:\cdots : s_n]$.
Writing ${\bf x}$ and ${\bf s}$ for the corresponding column vectors of length $n+1$,
we represent 
\[\,\L_{A}^\perp +{\bf g}\;\;\text{ by }\;\;\{{\bf x} \in \PP^n: (-{\bf b},A) \cdot {\bf x} = 0\} \quad \text{and} \quad
\,\L_{A}+{\bf c}\;\;\text{ by }\;\;\{{\bf s} \in \PP^n : (-B {\bf c},B)\cdot{\bf s} = 0\}.\
\]
The projective primal-dual central curve $\mathcal{C} $
is an irreducible curve in $\PP^n \times \PP^n$. 
Its bi-homogeneous prime ideal
in $K[x_0,x_1,\ldots,x_n,s_0,s_1,\ldots,s_n]$ can be computed by the process
of saturation. Namely, we compute it as
 the saturation with respect to $\langle x_0 s_0 \rangle$
of the ideal generated by the above linear forms together with the bi-homogeneous forms
$x_i s_i - x_j s_j$. The resulting ideal is irreducible because the vectors ${\bf c}$ and ${\bf g}$ are generic.

\begin{example}
\label{ex:saturate}
Let $d=2,n=4$. Fix $2 {\times} 4$-matrices $A = (a_{ij})$ and $B = (b_{ij})$ such that
$\L_B=\L_A^{\perp}$. 
We start with  the ideal $J$ in $K[x_0,\ldots,x_4,s_0, \ldots,s_4]$ generated by
$$
\begin{matrix}
a_{11} (x_1-g_1 x_0) + a_{12} (x_2-g_2 x_0) + a_{13} (x_3-g_3 x_0) + a_{14}( x_4 -g_4 x_0), \\
a_{21} (x_1-g_1 x_0) + a_{22} (x_2-g_2 x_0) + a_{23} (x_3-g_3 x_0) + a_{24}( x_4 -g_4 x_0) ,\\
b_{11} (s_1-c_1 s_0) + b_{12} (s_2-c_2 s_0) + b_{13} (s_3-c_3 s_0) + b_{14}( s_4 -c_4 s_0) ,\\
b_{21} (s_1-g_1 s_0) + b_{22} (s_2-c_2 s_0) + b_{23} (s_3-c_3 s_0) + b_{24}( s_4 -c_4 s_0) ,\\
s_1 x_1 - s_2 x_2 , \,\,\;\;\;
s_2 x_2 - s_3 x_3 , \,\, \;\;\;
s_3 x_3 - s_4 x_4 .
\end{matrix}
$$
The central curve $\mathcal{C}$ is irreducible  in $\PP^4 \times \PP^4$.
It has degree $(3,3)$ unless $A$ is very special.
The prime ideal of $\mathcal{C}$ is computed as the saturation $\,(J : \langle x_0 s_0 \rangle^\infty) $.
We find that this ideal 
has two minimal generators in addition to the seven above.
These are cubic polynomials in ${\bf x}$ and in
${\bf s}$, which define the primal and dual central curves.
They are shown
in Figure~\ref{fig:global}.
\hfill $\diamond$
\end{example}

Returning to the general case, we note that, since our curve $\mathcal{C}$ lives in $\PP^n \times \PP^n$,  
its {\em degree} is now a pair of integers $(d_{\bf x}, d_{\bf s})$.
These two integers can be defined geometrically:
$d_{\bf x}$ is the number of solutions of a 
general equation $  \sum_{i=0}^n \alpha_i x_i = 0 $ on the curve $\mathcal{C}$,
and $d_{\bf s}$ is the number of solutions 
 of a general equation $  \sum_{i=0}^n \beta_i s_i = 0 $  on $\mathcal{C}$.
 
\begin{corollary} \label{prop:analycenter}
Let ${\bf c}$ and ${\bf g}$ be generic vectors in $\R^n$
and let $(d_{\bf x}, d_{\bf s})$ be the degree of the projective
 primal-dual central curve $\mathcal{C} \subset \PP^n \times \PP^n$.
 This degree is given by our two M\"obius numbers, namely
 $\,d_{\bf x} = |\mu(A,{\bf c})|$ and  $d_{\bf s} = |\mu(B,{\bf g})|$.
 These numbers are defined in {\rm (\ref{eq:moebius6})}.
\end{corollary}

\begin{proof}
The projection from the primal-dual central curve onto its image in either
${\bf x}$-space or ${\bf s}$-space is birational. For instance, if ${\bf x}$ is
a general point on the primal central curve then the corresponding point ${\bf s}$
is uniquely obtained by solving the linear equations $x_i s_i = x_j s_j$
on $\L_{A}+{\bf c}$. Likewise, given a general point ${\bf s}$ on the dual central curve
we can recover the unique ${\bf x}$ such that $({\bf x},{\bf s}) \in \mathcal{C}$.
This implies that the intersections in $\PP^n \times \PP^n$ that define $d_{\bf x}$ and $d_{\bf s}$ 
are equivalent to intersecting
the primal or dual central curve with a general hyperplane in $\PP^n$, and the number
of points on that intersection is the respective M\"obius number.
\end{proof}

Next we discuss the geometry of this correspondence between  the primal and dual curves 
at their special points, namely vertices and analytic centers of the relevant hyperplane arrangements.
These special points are given by intersecting the primal-dual curve $\C$ with certain bilinear equations. 
%$\{x_1s_1=0\}$ and  $\{x_0s_0=0\}$, respectively.
The sum of our two M\"obius numbers, $d_{\bf x}+d_{\bf s}$, is the
number of solutions of a general bilinear equation $\,\sum_{i,j} \gamma_{ij} x_i s_j = 0\,$
on the primal-dual central curve $\mathcal{C}$. Two special choices of
such bilinear equations are of particular interest, namely, the
bilinear equation $\,x_0 s_0 = 0 \,$ and the bilinear equation $\,x_i s_i = 0 \,$
for any $i \geq 1$. Note that the choice of the index $i$ does not matter for the
second equation because $x_i s_i = x_j s_j$ holds on the curve.

Let us first observe what happens in $\PP^n \times \PP^n$
when the parameter $\lambda$ becomes $0$.
The corresponding points on the primal-dual curve $\mathcal{C}$
are found by solving the equation $x_1 s_1  = 0$ on
$\mathcal{C}$.  Its points
are the solutions of the $n$ equations
 $\,x_1 s_1 = x_2 s_2 = \cdots = x_n s_n = 0\,$
on  the $n$-dimensional subvariety
$\,  (\L_{A}^\perp + {\bf g}) \times (\L_{A} + {\bf c})\,$ of $\,\PP^n \times \PP^n$.
This intersection now contains many points in the product of
affine spaces, away from the hyperplanes $\{x_0 = 0\}$ and $\{s_0 = 0\}$.
We find the points by solving the linear equations
$x_{i_1}  = \cdots = x_{i_d} = 0$ on $\L_{A}^\perp + {\bf g}$ and 
$s_{j_1} = \cdots = s_{j_{n-d}} = 0$ on $\L_{A} + {\bf c}$,
where $\{i_1,\ldots,i_d\}$ runs over all bases
of the matroid $M(\L_A)$ and $\{j_1,\ldots,j_{n-d}\}$ is the complementary
basis of the dual matroid $M(\L_A)^* = M(\L_B)$.
These points represent vertices in the hyperplane arrangements $\mathcal{H}$
and $\mathcal{H}^*$, where
\[\mathcal{H} \;\text{ denotes } \;\{x_i=0\}_{i\in[n]} \,\text{ in }\, \L_{A}^\perp + {\bf g} \;\;\;\;\text{ and }\;\;\;\;
\mathcal{H^*}\; \text{ denotes }\; \{s_i=0\}_{i\in[n]} \,\text{ in }\, \L_{A} + {\bf c}.\] 
The vertices come in pairs corresponding to complementary bases, so
  the points with parameter $\lambda = 0$ on the primal-dual central curve $\mathcal{C}$
are the pairs $({\bf x},{\bf s})$ where
${\bf x}$ is a vertex in the hyperplane arrangement $\mathcal{H}$
and ${\bf s}$ is the complementary vertex in the dual arrangement~$\mathcal{H}^*$.

Imposing the equation $x_0 s_0 = 0$ means setting $\lambda = \infty$ in the
parametric representation of the central curve. The points thus
obtained have the following geometric interpretation 
in terms of bounded regions of the hyperplane arrangements $\mathcal{H}$ and $\mathcal{H}^*$.
We recall that the \textit{analytic center} of the polytope $P=\{A \bf x=\bf b, \bf x \geq 0\}$
 is the unique point
in the interior of $P$ that maximizes the concave function $\sum_{i=1}^n \log(x_i)$.
The algebraic characterization of the analytic center is that the 
gradient of $\sum_{i=1}^n \log(x_i)$, which is ${\bf x}^{-1}$, is orthogonal to the affine space 
$\L_{A}^\perp + {\bf g} = \{A {\bf x}={\bf b}\}$. This means that the vector
${\bf x}^{-1}$  lies in the row span $\L_{A}$ of $A$. Let $\L_{A}^{-1}$ denote  the
   the coordinate-wise reciprocal of $\L_{A}$, regarded as a subvariety  in $\PP^n$.

\begin{proposition} \label{cor:ds}
The intersection $\,\L_{A}^{-1} \,\cap \, (\L_{A}^\perp + {\bf g}) \,$ defines
 a zero-dimensional reduced subscheme of the affine space $\PP^n \backslash \{x_0 = 0\}$. 
All its points are defined over $\R$. They are the analytic centers of the
polytopes that form the bounded regions of the arrangement~$\mathcal{H}$.
\end{proposition}

\begin{proof}
The analytic center of each bounded region is a point in the
variety  $\,\L_A^{-1} \,\cap \, (\L_A^\perp + {\bf g}) $, by the gradient 
argument in the paragraph above.
This gives us $|\mu(A)|$ real points of intersection on $\,\L_{A}^{-1} \,\cap \, (\L_{A}^\perp + {\bf g}) \,$.
By replacing $\L_{A,{\bf c}}^{-1}$ with $\L_{A}^{-1}$ in Proposition~\ref{prop:PS}, we know that the degree 
of $\L_A^{-1}$ is $|\mu(A)|$. This shows that these
real points are all the intersection points
(over $\mathbb{C}$) and they occur with multiplicity one.
This argument closely follows the proof of Lemma~\ref{lem:zwei}.
\end{proof}

\begin{figure}
\includegraphics[width=7.4cm]{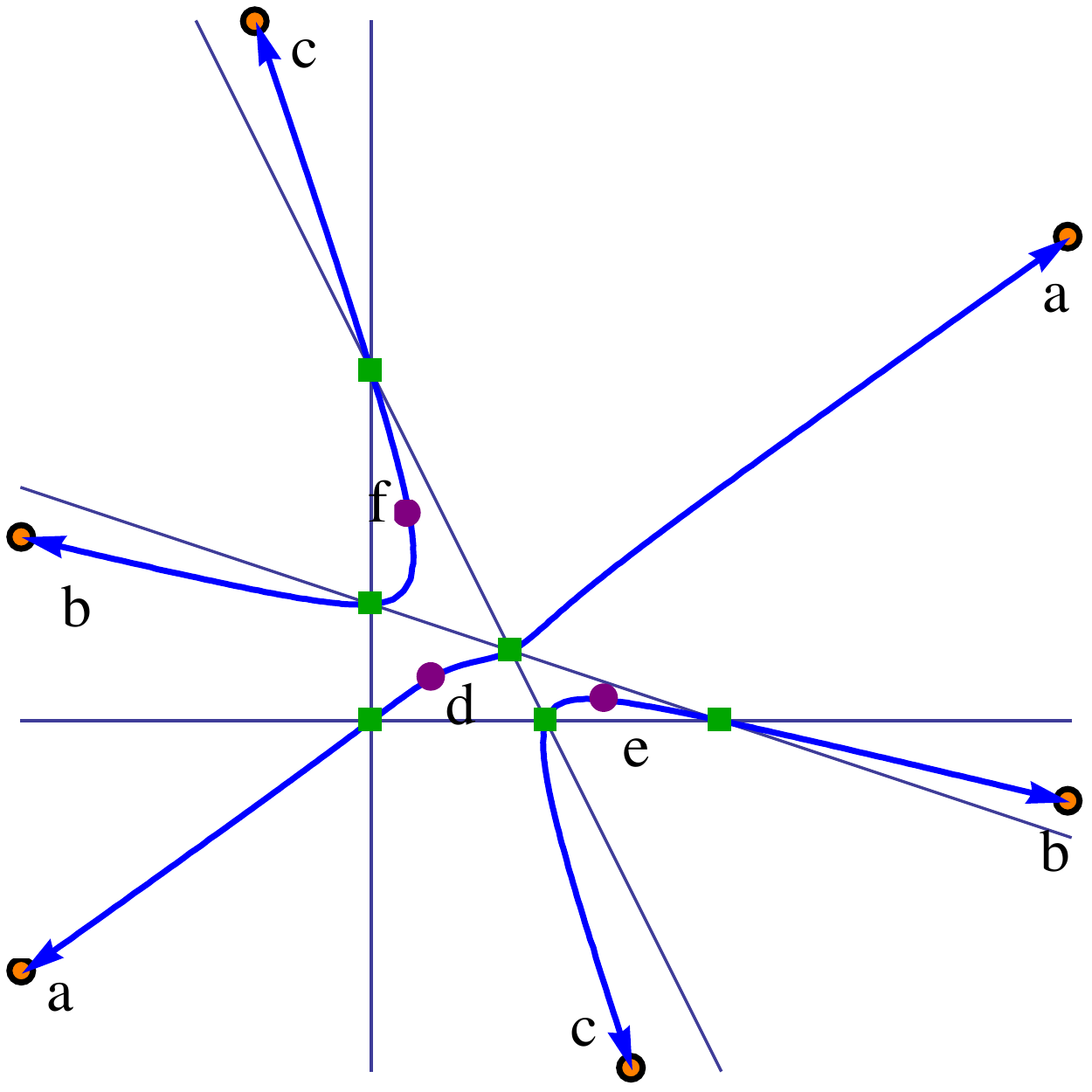}\hspace{0.7cm}
\includegraphics[width=7.4cm]{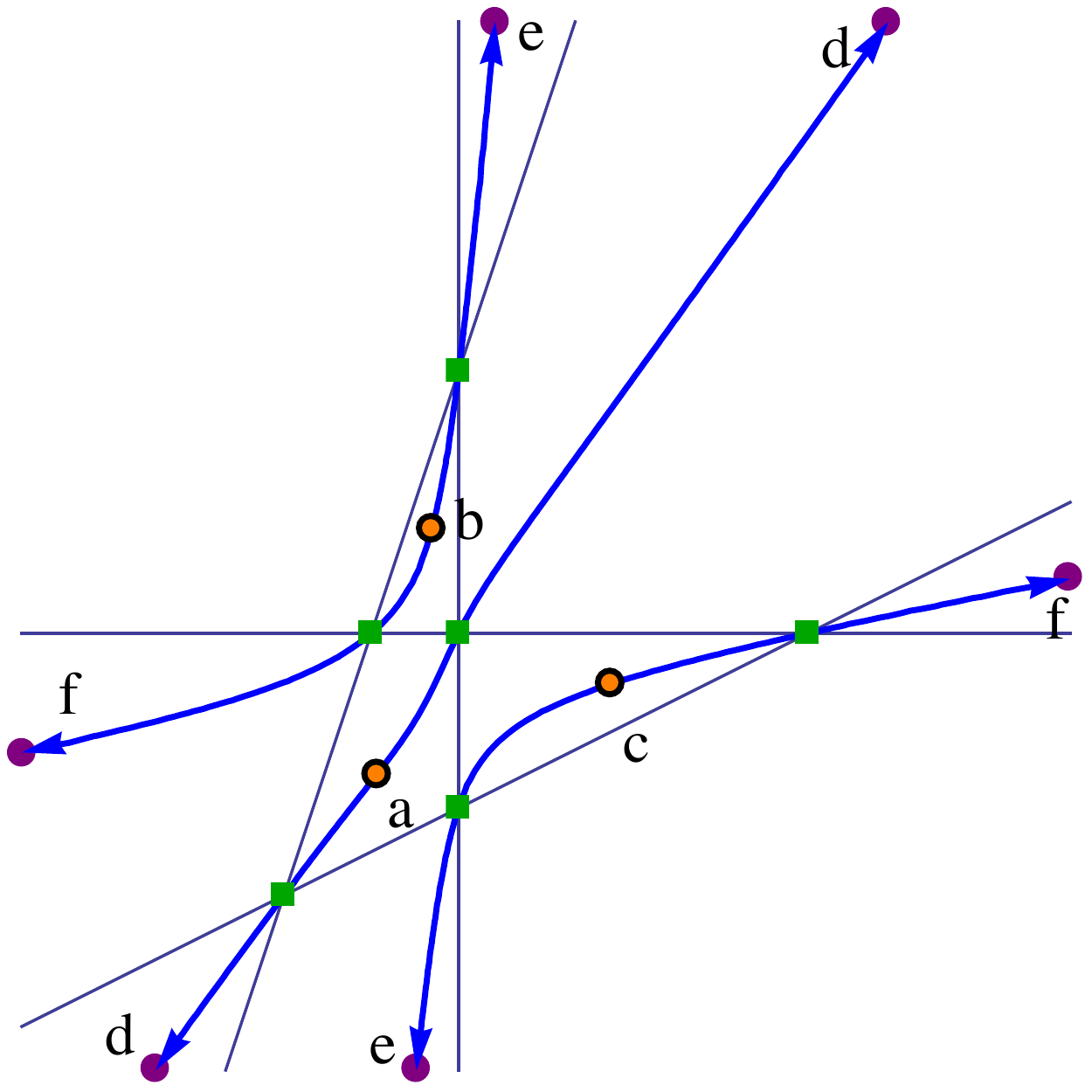}
\caption{Correspondence of vertices and analytic centers in the two projections of a primal-dual central curve.
Here both curves are plane cubics.}
\label{fig:global}
\end{figure}

Naturally, the dual statement holds verbatim, and we shall now state it explicitly.

\begin{proposition} \label{cor:dx}
The intersection $\,(\L_A^\perp)^{-1} \,\cap \, (\L_A + {\bf  c}) \,$ 
defines a zero-dimensional reduced subscheme of the affine space $\PP^n \backslash \{s_0 = 0\}$. 
All its points are  defined over $\R$.  They are the analytic centers of the
polytopes that form the bounded regions of the dual arrangement~$\mathcal{H}^*$.
\end{proposition}

We summarize our discussion with the following theorem on the
global geometry of the primal-dual central curve.
Figure \ref{fig:global} serves as an illustration for the case $n=4$ and $d=2$.

\begin{theorem}\label{thm:global}
The primal central curve in ${\bf x}$-space passes through each
vertex of the arrangement $\mathcal{H}$ as the dual central curve in ${\bf s}$-space 
passes through the corresponding vertex of $\mathcal{H^*}$.
As the primal curve passes through the analytic center of each bounded region in $\mathcal{H}$, 
the dual curve reaches the hyperplane $\{s_0=0\}$.
Similarly, as the dual curve reaches the analytic center of each bounded region in $\mathcal{H^*}$, 
the primal curve meets the hyperplane $\{x_0=0\}$.
\end{theorem} 

\begin{remark} 
The picture painted by
Theorem \ref{thm:global} resembles the  unpublished results of Adler and Haimovich
 \cite{adler, borgwardt, haimovich} on {\em co-optimal paths}.
 They considered an LP with two objective functions ${\bf c}$ and  ${\bf c}'$,
 and they studied the parametric objective ${\bf c}+\rho {\bf c}'$. 
As $\rho$ runs from $-\infty$ to $\infty$, the optimal solutions form a path of edges
in the arrangement $\mathcal{H}$.
Adler and Haimovich showed that the average length (in edges) of this co-optimal path in a cell of
$\mathcal{H}$, conditional on the path being nonempty in that cell, is at most $O({\rm min}(d,n-d))$.
This bound is close to the  curvature bounds in \cite{DMS}.
It would be interesting to explore possible connections between 
co-optimal edge paths and the degenerations of central curves
constructed in Theorem \ref{thm:degen}.
\hfill $\diamond$
\end{remark}

The primal central curve misses precisely one of the antipodal pairs of unbounded regions of $\mathcal{H}$.
It corresponds to the region in the induced arrangement
at infinity that contains the point representing the cost function ${\bf c}$.
For a visualization see the picture of the central curve in Figure \ref{fig:polygons}.
Here a projective transformation of $\PP^2$ moves the line from infinity into $\R^2$.

The points described in Propositions \ref{cor:ds} and \ref{cor:dx} are precisely
those points on the primal-dual central curve $\mathcal{C}$ for which the parameter $\lambda$
becomes $ \infty$. Equivalently, in its embedding in $\PP^n \times \PP^n$,
these are  solutions of the equation $x_0 s_0 = 0$ on the curve $\mathcal{C}$.
Note, however, that for special choices of $A$, the projective curve $\mathcal{C}$
will pass though points with $x_0 = s_0 = 0$.
Such points, which lie on the hyperplanes at infinity in both projective spaces,
are entirely independent of the choice of ${\bf c}$ and ${\bf g}$.
Indeed, they are the solutions of the equations
\begin{equation}
\label{eq:disjointsupport}
  {\bf s} \in \L_{A} =\ker{B} \,,\,\,\;\;\;    {\bf x} \in \L_{A}^{\perp}=\ker{A}, \;\;
   \,\,\,  \,\text{and} \,\,\, \;\;\; x_1 s_1 = x_2 s_2 = \dots = x_n s_n = 0.
   \end{equation}
The solutions to these equations form the {\em disjoint support variety}
in $\PP^{n-1} \times \PP^{n-1}$, which contains
pairs of vectors in the two spaces $\L_{A}$ and $\L_{A}^\perp$ whose respective supports are disjoint.

\begin{figure}\begin{center}
\includegraphics[width=7cm]{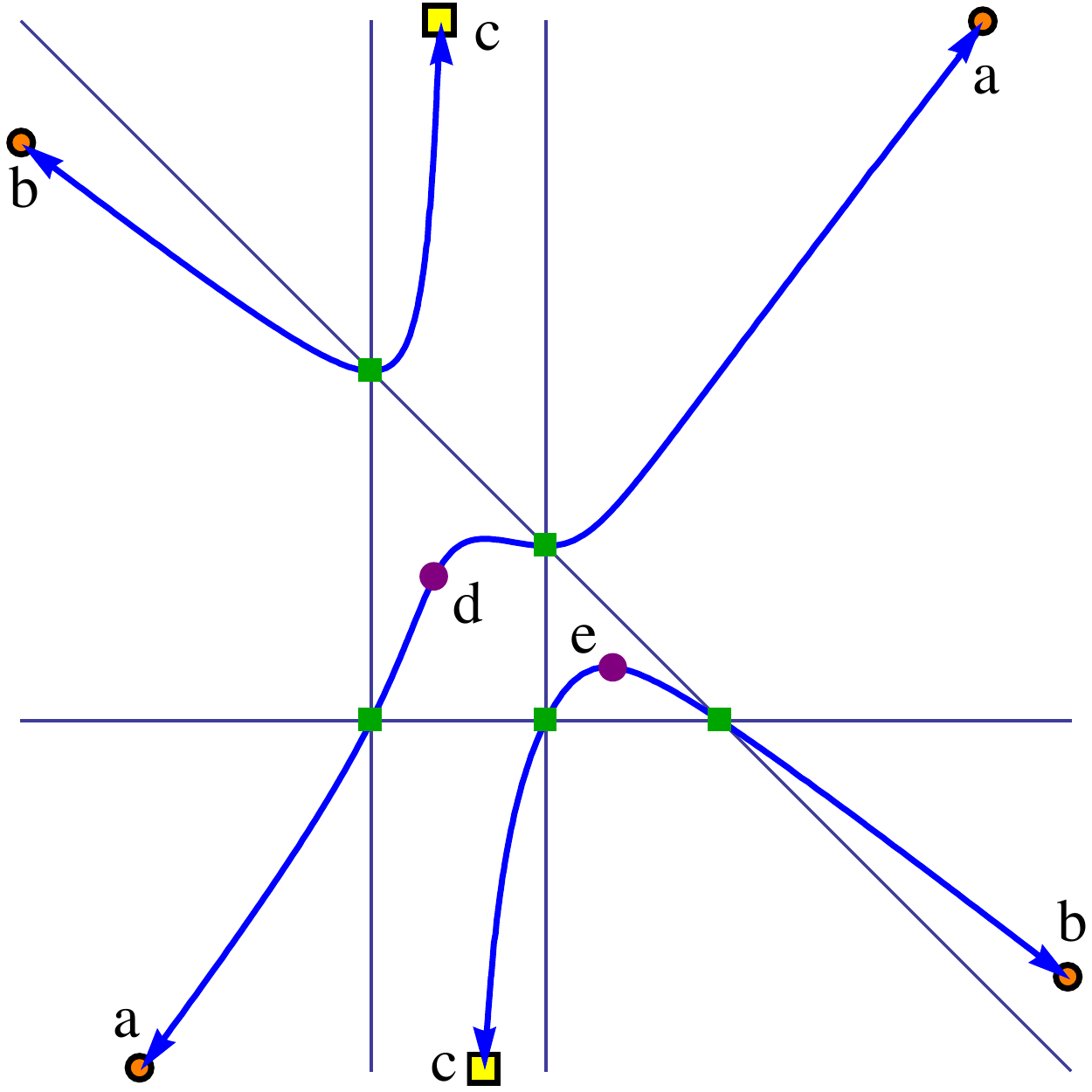}\hspace{0.7cm}
\includegraphics[width=7cm]{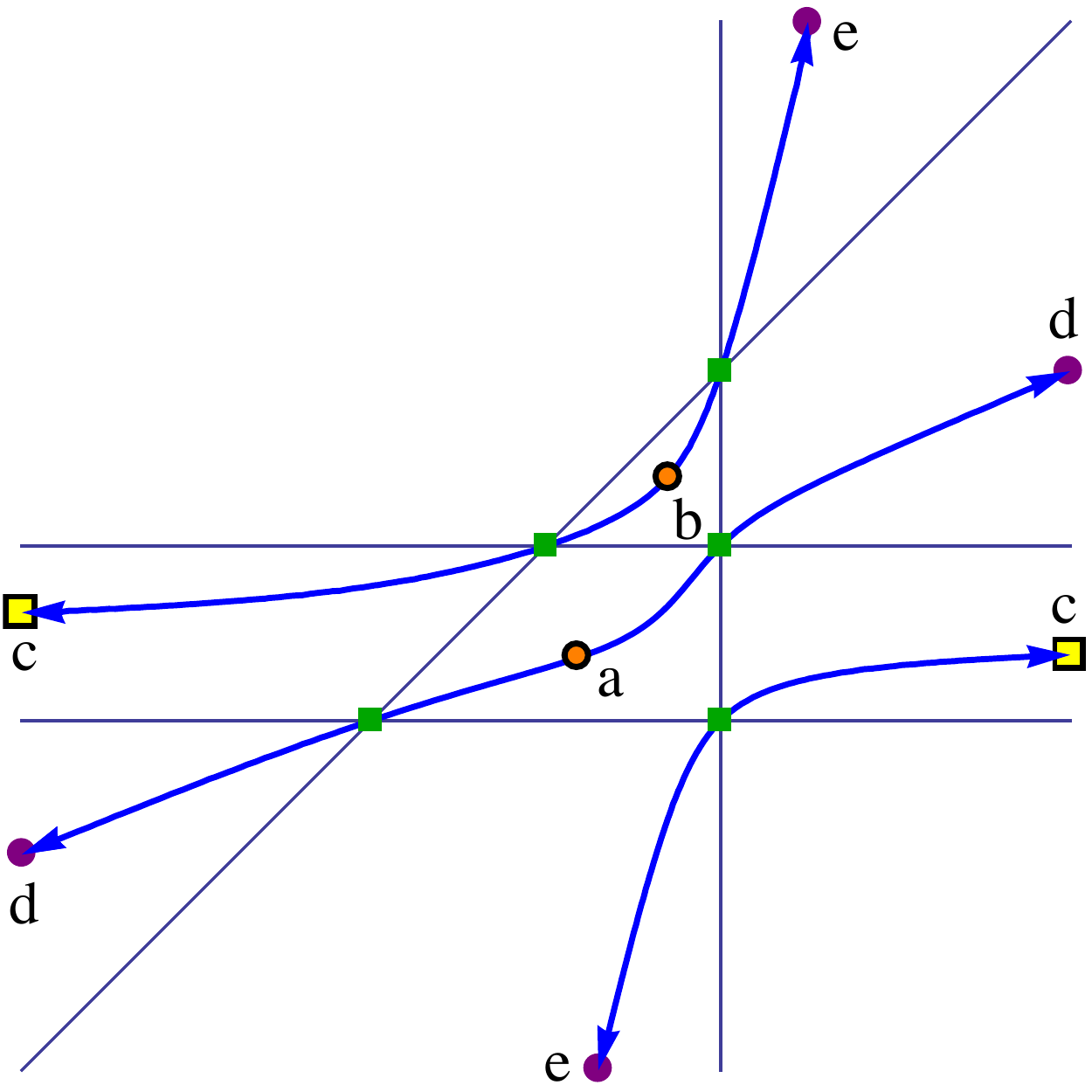} \end{center}
\caption{A primal-dual curve that
 intersects its disjoint support variety.}
\label{fig:dsv}
\end{figure}

\begin{example}\label{ex:dsv} Figure~\ref{fig:dsv} shows the primal-dual central curve
for the matrices
\[A = \begin{pmatrix}   1&-1&0&0\\ 0&1&1&-1 \end{pmatrix} 
\;\;\; \text{ and } \;\;\;
B = \begin{pmatrix} 1&1&0&1 \\ 0&0&1&1 \end{pmatrix}. \]
The disjoint support variety consists of the point 
$({\bf x}, {\bf s})= ([0:0:1:1], [1:-1:0:0])$  in $\PP^3\times \PP^3$. After we identify $\PP^3\times \PP^3$ with $\{({\bf x},{\bf s})\in \PP^4\times \PP^4 \;:\; x_0=0,s_0=0\}$, this point lies on the primal-dual central curve and appears as point ``c'' in Figure~\ref{fig:dsv}. 
\hfill $\diamond$
\end{example}

When studying the global geometry of the primal-dual central curve, it is useful
to start with the case when the constraint matrix $A$ is generic.  In that case, our
matroids are uniform, namely $M(\L_A) = U_{d,n}$ and $M(\L_B) = U_{n-d,n}$, and the
disjoint support variety (\ref{eq:disjointsupport}) is empty.
This condition ensures that the intersections of the curve $\mathcal{C}$
with both the hypersurfaces $\{x_0 s_0 = 0\}$ and $\{x_1 s_1 = 0\}$
in $\PP^n \times \PP^n$ is reduced, zero-dimensional and fully real.
The number of points in these intersections is the common number of bases in the two matroids:
$$ d_{\bf x} + d_{\bf s} \,\, = \,\, \binom{n-1}{d} + \binom{n-1}{d-1} \,\,\, = \,\,\,
\binom{n}{d} \,\, = \,\, \binom{n}{n-d} . $$
The intersection  points of $\mathcal{C}$ with $\{x_0 s_0 = 0\}$
are the pairs $({\bf x},{\bf s})$ where either ${\bf x}$ is an analytic center in $\mathcal{H}$
and ${\bf s}$ lies at infinity in the dual central curve, or
${\bf x}$ lies at infinity in the primal central curve and
${\bf s}$ is an analytic center in $\mathcal{H}^*$.
The intersection  points of $\mathcal{C}$ with $\{x_1 s_1 = 0\}$
are the pairs $({\bf x},{\bf s})$ where ${\bf x}$ is a vertex in $\mathcal{H}$
and ${\bf s}$ is a vertex in $\mathcal{H}^*$. 
Figure~\ref{fig:global} visualizes the above correspondences
for the case $n=4$ and $d=2$.
If  we now degenerate the generic matrix $A$
into a more special matrix, then some of the above points
representing vertices and analytic centers degenerate to points on the 
disjoint support variety (\ref{eq:disjointsupport}), as in Example~\ref{ex:dsv}.

In Theorem \ref{thm:global} we did not mention
the degree of the primal or dual central curve. For the sake of completeness,
here is a brief discussion of the geometric meaning of the degree~$d_{\bf x}$:

\begin{remark}\label{rem:CPtoAC}
Consider the intersection of the primal central path
with a  level set $\{{\bf c}^T {\bf x} =c_0\}$ of the 
objective function ${\bf c}$. Varying $c_0$ produces
a family of parallel hyperplanes. Each hyperplane
meets the curve in precisely $d_{\bf x}$ points,
all of which have real coordinates. These points
are the analytic centers of the $(n{-}d{-}1)$-dimensional 
polytopes obtained as the bounded regions of the induced
arrangement of $n$ hyperplanes $\{x_i = 0\}$ in the affine space
$\{{\bf x} \in \R^n :  A {\bf x} = {\bf b} ,\,
{\bf c}^T {\bf x} =c_0\}$. We can see $d_{\bf x}$ as the 
number of $(n-d-1)$-dimensional bounded regions in
the restriction of  the arrangement $\,\mathcal{H} \,$ to a general level hyperplane
$\{{\bf c}^T {\bf x} =c_0\}$. In particular, this gives a one-dimensional family of hyperplanes 
all of whose intersection points with the central curve are real,
as suggested by the left diagram in Figure \ref{fig:polygons}.
\hfill $\diamond$ 
\end{remark}

A main theme in this paper was that projective algebraic geometry provides 
an alternative view on optimality and duality in optimization, as well as powerful tools for analyzing 
interior point methods. This parallels the discussion of
semidefinite programming in \cite{philipbernd}. See also \cite{SU} for a
statistical perspective on analytic centers and central curves in the semidefinite context.

 \section{Acknowledgments}
We are grateful to Andrew Berget, Antoine Deza, Ragni Piene, Kristian Ranestad, Franz Rendl,
Raman Sanyal, Kim Chuan Toh, Michael J. Todd, Yuriy Zinchenko, and the referees for many  useful ideas, suggestions, 
clarifications,  and references.  Bernd Sturmfels and  Cynthia Vinzant acknowledge support by the 
U.S.~National Science Foundation (DMS-0757207 and DMS-0968882). 
Jes\'us De Loera was supported by the NSF grant DMS-0914107 and a UC Davis Chancellor Fellow award. 
He is also grateful to IPAM UCLA, where he was based during the writing of this article.

\bigskip
\end{document}